\documentclass[a4,numbers,sort&compress]{elsarticle}%
\usepackage{amsfonts}
\usepackage{amssymb,amsmath,amsthm,graphics,color}
\usepackage{graphics}
\usepackage{graphicx}
\usepackage{epsfig}
\usepackage{fancybox}
\usepackage{caption}
\usepackage{tikz,mathpazo}
\usepackage{multirow,booktabs}
\usepackage{geometry}
\usepackage{tabularx}
\usepackage{float}
\usepackage[justification=centering]{caption}
\usepackage[pagewise]{lineno}
\usepackage{amsmath}
\usepackage{amssymb}%
\usepackage{epstopdf}
\usetikzlibrary{shapes.geometric, arrows}
\usetikzlibrary{calc}
\journal{}

\newtheorem{theorem}{Theorem}[section]
\newtheorem{corollary}{Corollary}[section]

\newtheorem{lemma}{Lemma}[section]
\allowdisplaybreaks[1]
\newenvironment{lemmabis}[1]
{\renewcommand{\thelemma}{\ref{#1}$'$}
\addtocounter{lemma}{-1}
\begin{lemma}}
{\end{lemma}}
\begin{document}
\begin{frontmatter}
\title{Dynamics of COVID-19 models with asymptomatic infections and quarantine measures}
\author[mymainaddress,mysecondaryaddress]{Songbai Guo}
\ead{guosongbai@bucea.edu.cn}
\author[mymainaddress]{Yuling Xue}
\ead{xyl981274902@163.com}
\author[mysecondaryaddress1]{Xiliang Li\corref{mycorrespondingauthor}}
\ead{lixiliang@amss.ac.cn}
\author[mymainaddress0,mysecondaryaddress,mysecondaryaddress0]{Zuohuan Zheng}
\ead{zhzheng@amt.ac.cn}
\cortext[mycorrespondingauthor]{Corresponding author.}
\address[mymainaddress]{School of Science, Beijing University of Civil Engineering and Architecture, Beijing 102616, P. R. China}
\address[mymainaddress0]{School of Mathematics and Statistics, Hainan Normal University, Haikou 571158, P. R. China}
\address[mysecondaryaddress]{Academy of Mathematics and Systems Science, Chinese Academy of Sciences,
Beijing 100190, P. R. China}
\address[mysecondaryaddress1]{School of Mathematics and Information Science, Shandong Technology and Business University, Yantai 264005, P. R. China}
\address[mysecondaryaddress0]{School of Mathematical Sciences, University of Chinese Academy of Sciences,
Beijing 100049, P. R. China}
\begin{abstract}
Considering the propagation characteristics of COVID-19 in different regions, the dynamics analysis and numerical demonstration of long-term and short-term models of COVID-19 are carried out, respectively. The long-term model is devoted to investigate the global stability of COVID-19 model with asymptomatic infections and quarantine measures. By using the limit system of the model and Lyapunov function method, it is shown that the COVID-19-free equilibrium $V^0$ is globally asymptotically stable if the control reproduction number $\mathcal{R}_{c}<1$ and globally attractive if $\mathcal{R}_{c}=1$, which means that COVID-19 will die out; the COVID-19 equilibrium $V^{\ast}$ is globally asymptotically stable if $\mathcal{R}_{c}>1$, which means that COVID-19 will be persistent. In particular, to obtain the local stability of $V^{\ast}$, we use proof by contradiction and the properties of complex modulus with some novel details, and we prove the weak persistence of the system to obtain the global attractivity of $V^{\ast}$.
Moreover, the final size of the corresponding short-term model is calculated and the stability
of its multiple equilibria is analyzed. Numerical simulations of COVID-19 cases show that quarantine measures and asymptomatic infections have a non-negligible impact on the transmission of COVID-19.
\end{abstract}
\end{frontmatter}


\section{Introduction}

The newly discovered coronavirus disease 2019 (COVID-19) is a single-stranded
RNA coronavirus that can infect animals or human beings \cite{Xu20}. Recently, the
World Health Organization reported on 31 October 2022 that the cumulative
confirmed cases of COVID-19 in the world had exceeded 627 million, among which
more than 6.5 million had died from the virus \cite{WHO22}. The COVID-19 pandemic
affected the global economy and increased the economic burden on low-income
countries \cite{Nia22}. Following the discussion in \cite{Althobaity22} (also see \cite{Cui20,Liu21,Zhang22}), COVID-19 transmission occurred before the onset of symptoms. Research showed that asymptomatic infections could also cause COVID-19 transmission \cite{Gao07,Chen20}. The symptoms of COVID-19
transmission include cough, fever, fatigue, dyspnea and abdominal pain
\cite{Mizrahi20}. Susceptible individuals will be infected by contacting
infected individuals, inhaling virus-laden droplets, or touching the surface
of contaminated objects \cite{WHO20}. Public health and social measures, such
as wearing masks, disinfection, mass nucleic acid testing as well as quarantine measure, are
very helpful in controlling COVID-19 transmission \cite{WHO21}, and non-pharmaceutical interventions should not be relaxed prematurely \cite{Tang22}. Particularly, quarantine plays a special role in preventing further transmission of COVID-19 \cite{Basu22,Cui20,Giordano20,Roda20,Yuan20}. China takes a series of measures, including quarantine, which has made significant contributions to controlling the epidemic and reducing the fatality rate as stated in \cite{Ruan20}.

In order to make the media better disseminate information about COVID-19,
interdisciplinary methods including mathematics can be used to determine
appropriate communication strategies \cite{Briand21}. Moreover, the
establishment of appropriate COVID-19 mathematical models will help us to
understand the interplay between different pandemic factors \cite{Salman21}.
Kamara et al. \cite{Kamara21} considered a COVID-19 mathematical
model with the infectivity of exposed individuals, and analyzed the global
stability of disease-free and endemic equilibria, also see Bassey and Atsu
\cite{Bassey09}. Zamir et al. \cite{Zamir05} pointed out that exposed
individuals, symptomatic infected individuals, asymptomatic infected
individuals and stuffs contaminated with COVID-19 would infect susceptible
individuals. Cui et al. \cite{Cui20} presented a short-term model of COVID-19
transmission, and studied the final size of COVID-19 in Wuhan and
Guangzhou, respectively. Lv et al. \cite{Lv20} proposed long-term and
short-term mathematical models to illustrate the impact of asymptomatic
transmission on endemic. According to the research of McCallum et al.
\cite{McCallum01}, if each compartment in the model represented the actual
number of population rather than the density of population, the standard
incidence rate at this time would better represent the transmission rate of pathogens.

Recently, Guo et al. \cite{Guo22} established the following long-term COVID-19
model with quarantine and standard incidence rate,
\begin{equation}%
\begin{split}
\dot{S}(t)  &  =\lambda-\beta\frac{S(t)}{N(t)}(aE(t)+I(t)+bA(t))-dS(t),\\
\dot{E}(t)  &  =\beta\frac{S(t)}{N(t)}(aE(t)+I(t)+bA(t))-(c+d)E(t),\\
\dot{I}(t)  &  =pcE(t)-(q_{1}+r_{1}+d)I(t),\\
\dot{A}(t)  &  =(1-p)cE(t)-(q_{2}+r_{2}+d)A(t),\\
\dot{Q}(t)  &  =q_{1}I(t)+q_{2}A(t)-(r_{3}+d)Q(t),\\
\dot{R}(t)  &  =r_{1}I(t)+r_{2}A(t)+r_{3}Q(t)-dR(t),
\end{split}
\label{mod1}%
\end{equation}
where
\begin{equation}
N(t)=S(t)+E(t)+I(t)+A(t)+Q(t)+R(t), \label{nstt}%
\end{equation}
and the model parameters are all positive with definitions listed in the Tab.~\ref{tab1}.
They obtained the local asymptotic stability of the COVID-19-free equilibrium
$V^{0}=(S^{0},0,0,0,0,0)^{T}$~($S^{0}=\lambda/d$) of the model, and the
existence of COVID-19 equilibrium $V^{\ast}=(S^{\ast},E^{\ast},I^{\ast
},A^{\ast},Q^{\ast},R^{\ast})^{T}$ of the model in terms of control
reproduction number
\begin{equation}
\mathcal{R}_{c}=\frac{a\beta}{c+d}+\frac{pc\beta}{(c+d)B_{1}}+\frac
{bc\beta(1-p)}{(c+d)B_{2}}, \label{crn}%
\end{equation}
where $B_{i}=q_{i}+r_{i}+d$, $i=1,2$. In particular, they proposed a novel
analysis approach of uniform persistence of model \eqref{mod1} different from
the traditional persistence methods. It is not difficult to obtain that model
\eqref{mod1} is well-posed and dissipative in $D=\left\{
\varphi=(\varphi_{1},\varphi_{2},\varphi_{3},\varphi_{4},\varphi_{5}%
,\varphi_{6})^{T}\in\mathbb{R}_{+}^{6}:\sum_{i=1}^{6}\varphi_{i}>0\right\}  $
with $\mathbb{R}_{+}=[0,\infty)$ (see \cite{Guo22}).
\begin{table}[ptbh!]
\caption{Definitions of parameters in model \eqref{mod1}.}%
\label{tab1}%
\begin{tabular}
[c]{cl}%
\toprule Parameter & Definition\\
\midrule $\lambda$ & The birth rate of susceptible individuals\\
$d$ & The natural death rate\\
$\beta$ & The transmission rate of COVID-19\\
$a$ & The regulatory factor for infection probability of exposed individuals\\
$b$ & The regulatory factor for infection probability of asymptomatically
infected individuals\\
$c$ & The transfer rate of exposed individuals to other infected individuals\\
$p$ & The transition probability of symptomatically infected individuals\\
$q_{1}$ & The quarantined rate of symptomatically infected individuals\\
$q_{2}$ & The quarantined rate of asymptomatically infected individuals\\
$r_{1}$ & The recovery rate of symptomatically infected individuals\\
$r_{2}$ & The recovery rate of asymptomatically infected individuals\\
$r_{3}$ & The recovery rate of quarantined individuals\\\hline\hline
$S $ & Susceptible individuals\\
$E $ & Exposed individuals\\
$I $ & Symptomatically infected individuals\\
$A $ & Asymptomatically infected individuals\\
$Q $ & Quarantined individuals\\
$R $ & Recovered individuals\\
$N $ & Total population\\
\bottomrule &
\end{tabular}
\end{table}

Guo et al. \cite{Guo22} remarked that the global stability problems of $V^{0}$
and $V^{\ast}$ were very practical and challenging, and they would settle these
problems in future work. The purpose of the current research is to solve those
problems. Additionally, since the future development trend of COVID-19 is very
uncertain, the corresponding short-term model is established. We calculate the
control reproduction number and the final size of the model, and analyze the
stability of its multiple equilibria. Based on the characteristics of COVID-19
transmission in different regions, we choose long-term and short-term models for
numerical analysis. This method makes our research more practical and the
numerical fitting accuracy is higher. According to \cite{Lv20}, partial
derivatives of \eqref{crn} with respect to parameters can reflect some
information about COVID-19 transmission. Since this paper focuses on the impact of
asymptomatic infections $b$, $1-p$ and quarantine measures $q_{1}$, $q_{2}$ on the transmission of COVID-19, we calculate the partial
derivatives of equation \eqref{crn} with respect to the above four parameters as follows
\begin{align*}
\frac{\partial\mathcal{R}_{c}}{\partial b}  &  =\frac{c\beta(1-p)}{(c+d)B_{2}%
},~~ \frac{\partial\mathcal{R}_{c}}{\partial(1-p)} =-\frac{c\beta}{(c+d)B_{1}%
}+\frac{bc\beta}{(c+d)B_{2}},\\
\frac{\partial\mathcal{R}_{c}}{\partial q_{1}}  &  =-\frac{pc\beta}%
{(c+d)B_{1}^{2}},~~ \frac{\partial\mathcal{R}_{c}}{\partial q_{2}}
=-\frac{bc\beta(1-p)}{(c+d)B_{2}^{2}}.
\end{align*}
Obviously, the control reproduction number $\mathcal{R}_{c}$ is positively
correlated with parameter $b$, and negatively correlated with $q_{1}$ and
$q_{2}$. Furthermore, the relationship between $1-p$ and $\mathcal{R}_{c}$ is
related to the value range of $b$ under values of other parameters are fixed,
see Sections \ref{612} for specific details. We will elaborate on the
relationship among $\mathcal{R}_{c}$ and its parameters in Section \ref{LcsSA}
with a practical case.

The structure of this paper is as follows. In Section 2, the local stability
of COVID-19 equilibrium $V^{\ast}$ is proved by using proof by contradiction
and the properties of complex modulus. In Section 3, we obtain the weak
persistence of long-term COVID-19 model by using some analysis techniques. In
Section 4, the global stability of $V^{0}$ and $V^{\ast}$ is proved. In
Section 5, we propose a short-term COVID-19 model based on model \eqref{mod1}
and calculate its control reproduction number $\mathcal{R}_{c}$ and the final
size. At the same time, we analyze the stability of multiple equilibria of the
short-term model. In Section 6, the long-term and the short-term COVID-19
models are applied to case study in India and Nanjing, respectively, and the
sensitivity analysis of the corresponding $\mathcal{R}_{c}$ is carried out.
The impact of asymptomatic infections and quarantine measures on controlling
the spread of COVID-19 is discussed. The last section is the conclusions of
this paper.

\section{Local stability of the COVID-19 equilibrium}

From \cite[Theorem 4.1]{Guo22}, it follows that the COVID-19-free equilibrium
$V^{0}$ is locally asymptotically stable if the control reproduction number
$\mathcal{R}_{c}<1$ and unstable if $\mathcal{R}_{c}>1$. It is difficult to
prove the local stability of the COVID-19 equilibrium $V^{\ast}$ by using the
Routh-Hurwitz criterion. Hence, motivated by \cite{Alsh17}, our key idea to
settle this difficulty is to use the contradiction combining
the properties of complex modulus.

\begin{theorem}
\label{thm32} If $\mathcal{R}_{c}>1$, then the COVID-19 equilibrium $V^{\ast}$
is locally asymptotically stable.
\end{theorem}

\proof The characteristic equation at $V^{\ast}$ is given by
\begin{equation}
\left(  \Lambda+d\right)  \left(  \Lambda+r_{3}+d\right)  J=0, \label{eq321}%
\end{equation}
where%
\[%
\begin{split}
J=  &  \left[  \Lambda+\frac{d(\mathcal{R}_{c}-1)^{2}}{\mathcal{R}_{c}%
}+d\right]  \left[  \Lambda+(c+d)\right]  \left(  \Lambda+B_{1}\right)
\left(  \Lambda+B_{2}\right)  +pc\left(  \Lambda+d\right)  \left(
\Lambda+B_{2}\right)  \left[  \frac{d(\mathcal{R}_{c}-1)}{\mathcal{R}_{c}%
}-\frac{\beta}{\mathcal{R}_{c}}\right] \\
&  +\left(  \Lambda+d\right)  \left(  \Lambda+B_{1}\right)  \left(
\Lambda+B_{2}\right)  \left[  \frac{d(\mathcal{R}_{c}-1)}{\mathcal{R}_{c}%
}-\frac{\beta a}{\mathcal{R}_{c}}\right]  +(1-p)c\left(  \Lambda+d\right)
\left(  \Lambda+B_{1}\right)  \left[  \frac{d(\mathcal{R}_{c}-1)}%
{\mathcal{R}_{c}}-\frac{\beta b}{\mathcal{R}_{c}}\right]  .
\end{split}
\]
Obviously, equation \eqref{eq321} has negative roots $-d$ and $-(r_{3}+d)$.
The other eigenvalues satisfy $J=0$.

Next, we show that any root $\Lambda$ of $J=0$ has negative real part by
contradiction. Assume that $\Lambda$ has a non-negative real part. Then
$\frac{1}{\Lambda+y}$ also has a non-negative real part for $y\geq0$. Thus,
there holds
\begin{equation*}
\left[  \Lambda+(c+d)\right]  \left[ 1+\frac{d(\mathcal{R}_{c}-1)^{2}%
}{\mathcal{R}_{c}(\Lambda+d)}\right]  +\left[  \frac{pc}{\Lambda+B_{1}%
}+1+\frac{(1-p)c}{\Lambda+B_{2}}\right]  \frac{d(\mathcal{R}_{c}%
-1)}{\mathcal{R}_{c}}=\left[  \frac{pc\beta}{\Lambda+B_{1}}+\beta
a+\frac{(1-p)c\beta b}{\Lambda+B_{2}}\right]  \frac{1}{\mathcal{R}_{c}}.
\label{eq323}%
\end{equation*}
From the properties of complex modulus and $\mathcal{R}_{c}>1$, it follows that
\[
\left\vert \left[  \Lambda+(c+d)\right]  +\frac{d(\mathcal{R}_{c}-1)^{2}%
}{\mathcal{R}_{c}}+\frac{c}{\Lambda+d}\frac{d(\mathcal{R}_{c}-1)^{2}%
}{\mathcal{R}_{c}}+\left[  \frac{pc}{\Lambda+B_{1}}+1+\frac{(1-p)c}%
{\Lambda+B_{2}}\right]  \frac{d(\mathcal{R}_{c}-1)}{\mathcal{R}_{c}%
}\right\vert >c+d,
\]
and
\begin{align*}
&  \left\vert \left[  \frac{pc\beta}{\Lambda+B_{1}}+\beta a+\frac{(1-p)c\beta
b}{\Lambda+B_{2}}\right]  \frac{1}{\mathcal{R}_{c}}\right\vert \\
=  &  \left\vert \frac{pc\beta}{(c+d)(\Lambda+B_{1})}+\frac{\beta a}%
{c+d}+\frac{(1-p)c\beta b}{(c+d)(\Lambda+B_{2})}\right\vert \frac
{c+d}{\mathcal{R}_{c}}\\
\leq &  \left(  \left\vert \frac{pc\beta}{(c+d)(\Lambda+B_{1})}\right\vert
+\frac{\beta a}{c+d}+\left\vert \frac{(1-p)c\beta b}{(c+d)(\Lambda+B_{2}%
)}\right\vert \right)  \frac{c+d}{\mathcal{R}_{c}}\\
\leq &  \left[  \frac{pc\beta}{(c+d)B_{1}}+\frac{\beta a}{c+d}+\frac
{(1-p)c\beta b}{(c+d)B_{2}}\right]  \frac{c+d}{\mathcal{R}_{c}}=c+d.
\end{align*}
Clearly, this is a contradiction, and hence any root of equation \eqref{eq321}
has negative real part. Therefore, $V^{\ast}$ is locally asymptotically stable for $\mathcal{R}_{c}>1$.

\section{Weak persistence}

In this section, we will study the weak persistence of model \eqref{mod1}
based on some analysis techniques in \cite{Guo16,Guo22}. Let $\Omega
=\{\varphi\in\mathbb{R}_{+}^{6}:\varphi_{2}>0\}$ and
\[
u(t)\equiv(u_{1}(t),u_{2}(t),u_{3}(t),u_{4}(t),u_{5}(t),u_{6}(t))^{T}%
=(S(t),E(t),I(t),A(t),Q(t),R(t))^{T}%
\]
be the solution of model \eqref{mod1} with any $\varphi\in\Omega.$ It follows
that $\Omega\subseteq D$ is a positive invariant set for model \eqref{mod1},
and $u(t)\gg\mathbf{0}$ for $t>0$. Subsequently, we thus discuss the weak
persistence of model \eqref{mod1} in $\Omega$.

Model \eqref{mod1} is called weakly persistent if $\limsup_{t\rightarrow
\infty}u_{i}(t)>0$, $i=1,2,3,4,5,6$ for any $\varphi\in\Omega$ (see
\cite{Butler86}). To start the weak persistence of model \eqref{mod1}, the
following lemma is needed.

\begin{lemma}
\label{lem01} If $\mathcal{R}_{c}>1$, $\theta\in(0,1)$, and $\limsup
_{t\rightarrow\infty}E(t)\leq\theta E^{\ast}$, then it follows%
\[
\liminf_{t\rightarrow\infty}S(t)\geq\bar{S}\equiv\frac{\lambda}{\theta
\beta\left(  aE^{\ast}+I^{\ast}+bA^{\ast}\right)  /S^{0}+d}=\frac{S^{0}%
}{\theta\left(  \mathcal{R}_{c}-1\right)  +1}>S^{\ast}.
\]

\end{lemma}

\proof First, from the positive equilibrium equations, we can obtain $E^{\ast
}=\lambda(\mathcal{R}_{c}-1)/(c+d)\mathcal{R}_{c}.$ Thus, there holds $\bar
{S}=S^{0}/[\theta\left(  \mathcal{R}_{c}-1\right)  +1]>S^{\ast}.$ By the third
and the fourth equations of model \eqref{mod1}, we have that%
\[
\operatorname*{limsup}_{t\rightarrow\infty}I(t)\leq\theta I^{\ast},\text{
}\operatorname*{limsup}_{t\rightarrow\infty}\text{{}}A(t)\leq\theta A^{\ast}.
\]
It follows from model \eqref{mod1} that $\dot{N}(t)=\lambda-dN(t),$ which
yields
\begin{equation}
\lim_{t\rightarrow\infty}N(t)=S^{0}. \label{ntt}%
\end{equation}
Hence, we have
\[
\operatorname*{limsup}_{t\rightarrow\infty}\frac{aE(t)+I(t)+bA(t)}{N(t)}%
\leq\frac{\theta\left(  aE^{\ast}+I^{\ast}+bA^{\ast}\right)  }{S^{0}}\text{,}%
\]
Consequently, the first equation of model \eqref{mod1} implies
\[
\liminf_{t\rightarrow\infty}S(t)\geq\frac{\lambda}{\theta\beta\left(
aE^{\ast}+I^{\ast}+bA^{\ast}\right)  /S^{0}+d}=\bar{S}.
\]

In fact, the result of Lemma \ref{lem01} can be further improved. We have the
following Lemma.

\begin{lemmabis}{lem01}
\label{lem010} Under the conditions of Lemma \ref{lem01}, it holds%
\[
\liminf_{t\rightarrow\infty}S(t)\geq\breve{S}\equiv\frac{\lambda
-\theta(c+d)E^{\ast}}{d}=\frac{S^{0}}{1+\theta\left(  \mathcal{R}
_{c}-1\right)  /[\mathcal{R}_{c}-\theta\left(  \mathcal{R}_{c}-1\right)
]}>\bar{S}.
\]

\end{lemmabis}

\proof From the first two equations of model \eqref{mod1}, we have
\[
\liminf_{t\rightarrow\infty}(S(t)+E(t))\geq\frac{\lambda-c\theta E^{\ast}}%
{d}.
\]
In consequence,
\[
\liminf_{t\rightarrow\infty}S(t)\geq\liminf_{t\rightarrow\infty}%
(S(t)+E(t))-\limsup_{t\rightarrow\infty}E(t)\geq\breve{S},
\]
and the positive equilibrium equations imply%
\[
\breve{S}=\frac{S^{0}}{1+\theta\left(  \mathcal{R}_{c}-1\right)
/[\mathcal{R}_{c}-\theta\left(  \mathcal{R}_{c}-1\right)  ]}.
\]
Observe that $\mathcal{R}_{c}-\theta\left(  \mathcal{R}_{c}-1\right)  >1$, we
thus have $\breve{S}>\bar{S}.$

\begin{theorem}
\label{thm1} If $\mathcal{R}_{c}>1$ and $\theta\in(0,1)$, then $\limsup
_{t\rightarrow\infty}E(t)>\theta E^{\ast}$.
\end{theorem}

\proof We prove the result by contradiction. Suppose $\limsup_{t\rightarrow
\infty}E(t)\leq\theta E^{\ast}$. Then it follows from Lemma \ref{lem01} that
there is an $\varepsilon_{0}>0$ such that%
\begin{equation}
\frac{\bar{S}}{S^{0}+\varepsilon}>\frac{S^{\ast}}{S^{0}} \label{stt}%
\end{equation}
for any $\varepsilon\in(0,\varepsilon_{0})$. By Lemma \ref{lem01} and
\eqref{ntt}, we have that for any $\varepsilon\in(0,\varepsilon_{0}),$ there
exists a $T\equiv T\left(  \varepsilon,\varphi\right)  >0$ such that for all
$t\geq T$, there holds%
\[
\frac{S(t)}{N(t)}>\frac{\bar{S}}{S^{0}+\varepsilon}.
\]
Now, we define a function as follows
\[
L(\varphi)=\varphi_{2}+\frac{\beta\bar{S}}{B_{1}\left(  S^{0}+\varepsilon
\right)  }\varphi_{3}+\frac{b\beta\bar{S}}{B_{2}\left(  S^{0}+\varepsilon
\right)  }\varphi_{4},~\varphi\in\Omega\text{.}%
\]
Then the derivative of $V$ along the solution $u(t)$ for $t\geq T$ is given by%
\[
\dot{L}(u(t))\geq(c+d)\left(  \frac{\bar{S}}{S^{0}+\varepsilon}\mathcal{R}%
_{c}-1\right)  E(t)\text{.}%
\]

Denote
\[
m=\min\left\{  E(T),\frac{B_{1}I(T)}{pc},\frac{B_{2}A(T)}{(1-p)c}\right\}
\text{.}%
\]
Next, we will prove $E(t)\geq m$ for $t\geq T\text{.}$ If not, there exists a
$T_{0}\geq0$ such that $E(t)\geq m$ for $t\in\lbrack T,T+T_{0}]$,
$E(T+T_{0})=m$ and $\dot{E}(T+T_{0})\leq0\text{. }$For $t\in\lbrack
T,T+T_{0}]$, we have that
\begin{equation}
\dot{I}(t)=pcE(t)-B_{1}I(t)\geq pcm-B_{1}I(t). \label{etm}%
\end{equation}
It follows from \eqref{etm} that for $t\in\lbrack T,T+T_{0}]$,
\[
I(t)\geq\dfrac{pcm}{B_{1}}+\left(  I(T)-\dfrac{pcm}{B_{1}}\right)
e^{B_{1}(T-t)}\geq\dfrac{pcm}{B_{1}}.
\]
Similarly, we have
\[
A(t)\geq\frac{(1-p)cm}{B_{2}}%
\]
for $t\in\lbrack T,T+T_{0}]$. By Remark 3.1 in \cite{Guo22} and \eqref{stt},
we have
\[
\frac{\bar{S}}{S^{0}+\varepsilon}\mathcal{R}_{c}-1>\frac{S^{\ast}}{S^{0}%
}\mathcal{R}_{c}-1=0.
\]
In consequence, it holds that
\[
\dot{E}(T+T_{0})\geq m(c+d)\left(  \frac{\bar{S}}{S^{0}+\varepsilon
}\mathcal{R}_{c}-1\right)  >m(c+d)\left(  \frac{S^{\ast}}{S^{0}}%
\mathcal{R}_{c}-1\right)  =0.
\]
Clearly, this contradicts $\dot{E}(T+T_{0})\leq0$. Consequently, $E(t)\geq m$
for $t\geq T\text{.}$ Hence, it follows for $t\geq T$ that
\[
\dot{L}(u(t))\geq(c+d)\left(  \frac{\bar{S}}{S^{0}+\varepsilon}\mathcal{R}%
_{c}-1\right)  m>0,
\]
which hints $L(u(t))\rightarrow\infty$ as $t\rightarrow\infty\text{.}$
Therefore, this contradicts the boundedness of $L(u(t))$.

From Theorem \ref{thm1}, it is not difficult to obtain the following corollary.

\begin{corollary}
If $\mathcal{R}_{c}>1$, then model \eqref{mod1} is weakly persistent.
\end{corollary}

\section{Global stability}

In this section, we will discuss the global asymptotic stability of
COVID-19-free equilibrium $V^{0}$ and COVID-19 equilibrium $V^{\ast}$. Let
$u(t)\equiv(S(t),E(t),I(t),A(t),Q(t),R(t))^{T}$ be the solution of model
\eqref{mod1} with any $\varphi\in D$. Note that%
\[
\dot{N}(t)=\lambda-dN(t),
\]
we thus have $\lim_{t\rightarrow\infty}N(t)=S^{0}$. Then model \eqref{mod1}
has the following limiting system:%
\begin{equation}%
\begin{split}
\dot{S}_{1}(t)  &  =\lambda-\frac{\beta S_{1}(t)}{S^{0}}(aE_{1}(t)+I_{1}%
(t)+bA_{1}(t))-dS_{1}(t),\\
\dot{E}_{1}(t)  &  =\frac{\beta S_{1}(t)}{S^{0}}(aE_{1}(t)+I_{1}%
(t)+bA_{1}(t))-(c+d)E_{1}(t),\\
\dot{I}_{1}(t)  &  =pcE_{1}(t)-(q_{1}+r_{1}+d)I_{1}(t),\\
\dot{A}_{1}(t)  &  =(1-p)cE_{1}(t)-(q_{2}+r_{2}+d)A_{1}(t),\\
\dot{Q}_{1}(t)  &  =q_{1}I_{1}(t)+q_{2}A_{1}(t)-(r_{3}+d)Q_{1}(t),\\
\dot{R}_{1}(t)  &  =r_{1}I_{1}(t)+r_{2}A_{1}(t)+r_{3}Q_{1}(t)-dR_{1}(t).
\end{split}
\label{mod2}%
\end{equation}
It is not difficult to get that the solution $v(t)\equiv(S_{1}(t),E_{1}%
(t),I_{1}(t),A_{1}(t),Q_{1}(t),R_{1}(t))^{T}$ of model \eqref{mod2} through
any $\phi=(\phi_{1},\phi_{2},\phi_{3},\phi_{4},\phi_{5},\phi_{6})^{T}\in D$
exists, which is unique and nonnegative on $[0,\infty)$. Furthermore, we can get that
model \eqref{mod2} is dissipative in $D$ and $S_{1}(t)>0$ for $t>0$. Clearly,
$V^{0}$ and $V^{\ast}$ are also the equilibria of model \eqref{mod2}.

For the global stability of the COVID-19-free equilibrium $V^{0}$ of model \eqref{mod1},
we have the following conclusion.
\begin{theorem}
\label{thm51} The COVID-19-free equilibrium $V^{0}$ is globally asymptotically
stable if $\mathcal{R}_{c}<1$ and globally attractive if $\mathcal{R}_{c}=1$ in
$D$.
\end{theorem}

\proof We first know that $V^{0}$ is stable for $\mathcal{R}_{c}<1$ in the
light of \cite[Theorem 4.1]{Guo22}. Next, we will show that $V^{0}$ is
globally attractive for $\mathcal{R}_{c}\leq1$. Let $u(t)$ be the solution of
model \eqref{mod1} with any $\varphi\in D$. Since $u(t)$ is bounded, we can
get $\omega(\varphi)\subseteq D$ is compact, where $\omega(\varphi)$ is the
$\omega$-limit set of $\varphi$ with respect to model \eqref{mod1}. To prove that
$V^{0}$ is globally attractive, we only need to verify $\omega(\varphi
)=\{V^{0}\}.$

Let $v(t)$ be the solution of model \eqref{mod2} through any $\phi\in$ $D.$
Now, we use the technique of Lyapunov function in \cite[Theorem 4.3]{Bai21},
and then define the following function $L$ on $\mathcal{D}=\left\{  \phi
\in\mathbb{R}_{+}^{6}: \phi_{1}>0\right\}  \subseteq D$,%

\begin{equation}
L(\phi)=\phi_{1}-S^{0}-S^{0}\ln\frac{\phi_{1}}{S^{0}}+\phi_{2}+\frac{\beta
}{B_{1}}\phi_{3}+\frac{\beta b}{B_{2}}\phi_{4}. \label{vph}%
\end{equation}
It is easy to find that $L$ is continuous on $\mathcal{D}$. The derivative of $L$ along
this solution $v(t)$ for $t>0$ can be taken by
\begin{equation}%
\begin{split}
\dot{L}(v(t))  &  ={}\left(  1-\frac{S^{0}}{S_{1}(t)}\right)  \dot{S}%
_{1}(t)+\dot{E}_{1}(t)+\frac{\beta}{B_{1}}\dot{I}_{1}(t)+\frac{\beta b}{B_{2}%
}\dot{A}_{1}(t)\\
&  ={}\lambda-dS_{1}(t)-\frac{S^{0}\lambda}{S_{1}(t)}+dS^{0}+E_{1}(t)\left(
c+d\right)  \left(  \frac{a\beta}{c+d}+\frac{pc\beta}{B_{1}\left(  c+d\right)
}+\frac{b\beta\left(  1-p\right)  c}{B_{2}\left(  c+d\right)  }-1\right) \\
&  =-\frac{d}{S_{1}(t)}(S_{1}(t)-S^{0})^{2}+E_{1}(t)(c+d)(\mathcal{R}_{c}-1)\\
&  \leq{}0.
\end{split}
\label{vph1}%
\end{equation}
According to \eqref{vph} and \eqref{vph1}, $\omega(\phi)\subseteq D,$ where
$\omega(\phi)$ is the $\omega$-limit set of $\phi$ with respect to model
\eqref{mod2}. Thus, For $\mathcal{R}_{c}\leq1$, $L$ is a Lyapunov function on
$\{v(t):t\geq1\}\subseteq D$. It follows from \cite[Lemma 4.1]{Bai21} (also see
\cite[Corollary 2.1]{Guo19}) that $\dot{L}(\psi)=0$ for any $\psi\in
\omega(\phi)$.

For any $\psi\in\omega(\phi)$, let $v(t)$ be the solution of model
\eqref{mod2} through $\psi$. Then from the invariance of $\omega(\phi)$, it
follows that $v(t)\in\omega(\phi)$ for all $t\in\mathbb{R}$. It can be seen
from \eqref{vph1}, $S_{1}(t)=S^{0}$ for all $t\in\mathbb{R}$. By the first
equation of model \eqref{mod2}, we can obtain $E(t)=I(t)=A(t)=0$ for all
$t\in\mathbb{R}$. From the fifth and sixth equations of model \eqref{mod2} and
the invariance of $\omega(\phi)$, we have $Q(t)=R(t)=0$ for all $t\in
\mathbb{R}$. Consequently, $\omega(\phi)=\{V^{0}\}$ for $\mathcal{R}_{c}\leq
1$, and thus $W^{s}(V^{0})=D$, where $W^{s}(V^{0})$ is the stable set of
$V^{0}$ for model \eqref{mod2}. By the similar argument as in \cite[Theorem
4.3]{Bai21} , we can obtain $V^{0}$ is locally asymptotically stable for model
\eqref{mod2} if $\mathcal{R}_{c}\leq1$. It is clear to see that $\omega(\varphi)\cap W^{s}(V^{0}%
)\neq\emptyset.$ Therefore, \cite[Theorem 1.2]{Thieme92} implies that
$\omega(\varphi)=\{V^{0}\}.$

\begin{theorem}
\label{thm52} The COVID-19 equilibrium $V^{\ast}$ is globally asymptotically
stable if and only if $\mathcal{R}_{c}>1$ in $\Omega$.
\end{theorem}

\proof By Lemma 3.1 in \cite{Guo22}, we only need to verify its sufficiency. From
Theorem \ref{thm32}, $V^{\ast}$ is stable for $\mathcal{R}_{c}>1$. Let $u(t)$
be the solution of model \eqref{mod1} through any $\varphi\in\Omega$ and
$v(t)$ be the solution of model \eqref{mod2} through any $\phi\in\Omega$.
Notice that $\Omega$ is also positively invariant for model \eqref{mod2}, and $v(t)\gg\mathbf{0}$ for $t>0$. Let $\mathfrak{D}%
=\{\phi\in\mathbb{R}_{+}^{6}:\phi\gg0\}$. Then $\mathfrak{D\subseteq}\Omega.$
To show that $V^{\ast}$ is globally attractive for $\mathcal{R}_{c}>1$, we only
need to prove $\omega(\varphi)=\{V^{\ast}\}.$ Now, we use the technique of
Lyapunov function in \cite[Theorem 4.4]{Bai21} to define the following
function $L$ on $\mathfrak{D}$,%
\begin{equation}
L(\phi)=S^{\ast}g\left(  \frac{\phi_{1}}{S^{\ast}}\right)  +E^{\ast}g\left(
\frac{\phi_{2}}{E^{\ast}}\right)  +\frac{S^{\ast}\beta I^{\ast}}%
{S^{0}pcE^{\ast}}I^{\ast}g\left(  \frac{\phi_{3}}{I^{\ast}}\right)
+\frac{S^{\ast}\beta bA^{\ast}}{S^{0}(1-p)cE^{\ast}}A^{\ast}g\left(
\frac{\phi_{4}}{A^{\ast}}\right)  , \label{eq14}%
\end{equation}
where $g(x)=x-1-\ln x$, $x>0.$ The derivative of $L$ along $v(t)$ for $t>0$ is
given by
\[%
\begin{split}
\dot{L}(v(t))=  &  \left(  1-\frac{S^{\ast}}{S_{1}(t)}\right)  {}\dot{S}%
_{1}(t)+\left(  1-\frac{E^{\ast}}{E_{1}(t)}\right)  \dot{E}_{1}(t)+\frac
{S^{\ast}\beta I^{\ast}}{S^{0}pcE^{\ast}}\left(  1-\frac{I^{\ast}}{I_{1}%
(t)}\right)  \dot{I}_{1}(t)\\
&  +\frac{S^{\ast}\beta bA^{\ast}}{S^{0}(1-p)cE^{\ast}}\left(  1-\frac
{A^{\ast}}{A_{1}(t)}\right)  \dot{A}_{1}(t)\\
=  &  -\lambda\bigg[g\left(  \frac{S^{\ast}}{S_{1}(t)}\right)  +g\left(
\frac{S_{1}(t)}{S^{\ast}}\right)  \bigg]+(\frac{S^{\ast}\beta I^{\ast}}{S^{0}%
}+\frac{S^{\ast}\beta bA^{\ast}}{S^{0}}+\frac{S^{\ast}\beta aE^{\ast}}{S^{0}%
})g\left(  \frac{S_{1}(t)}{S^{\ast}}\right) \\
&  -\frac{S^{\ast}\beta I^{\ast}}{S^{0}}\left(  \frac{E_{1}(t)I^{\ast}%
}{E^{\ast}I_{1}(t)}+\frac{S_{1}(t)I_{1}(t)E^{\ast}}{S^{\ast}I^{\ast}E_{1}%
(t)}-\ln\frac{S_{1}(t)}{S^{\ast}}-2\right) \\
&  -\frac{S^{\ast}\beta bA^{\ast}}{S^{0}}\left(  \frac{E_{1}(t)A^{\ast}%
}{E^{\ast}A_{1}(t)}+\frac{S_{1}(t)A_{1}(t)E^{\ast}}{S^{\ast}A^{\ast}E_{1}%
(t)}-\ln\frac{S_{1}(t)}{S^{\ast}}-2\right) \\
&  -\frac{S^{\ast}\beta aE^{\ast}}{S^{0}}\left(  \frac{S_{1}(t)}{S^{\ast}%
}-1-\ln\frac{S_{1}(t)}{S^{\ast}}\right)  ,
\end{split}
\]
where the equilibrium equations of $V^{\ast}$:
\begin{align*}
d  &  =\frac{\lambda}{S^{\ast}}-\frac{\beta}{S^{0}}\left(  aE^{\ast}+I^{\ast
}+bA^{\ast}\right)  ,\\
c+d  &  =\frac{S^{\ast}\beta a}{S^{0}}+\frac{\beta S^{\ast}I^{\ast}}%
{S^{0}E^{\ast}}+\frac{\beta S^{\ast}bA^{\ast}}{S^{0}E^{\ast}},\\
q_{1}+d  &  =\frac{pcE^{\ast}}{I^{\ast}}-r_{1},\\
r_{2}+d  &  =(1-p)c\frac{E^{\ast}}{A^{\ast}}-q_{2}%
\end{align*}
are used. Further, we have%
\begin{equation}%
\begin{split}
\dot{L}(v(t)) =  &  -\lambda g\left(  \frac{S^{\ast}}{S_{1}(t)}\right)
-\left(  dS^{\ast}+\frac{S^{\ast}\beta aE^{\ast}}{S^{0}}\right)  g\left(
\frac{S_{1}(t)}{S^{\ast}}\right) \\
&  -\frac{S^{\ast}\beta I^{\ast}}{S^{0}}\left[  g\left(  \frac{E_{1}%
(t)I^{\ast}}{E^{\ast}I_{1}(t)}\right)  +g\left(  \frac{S_{1}(t)I_{1}%
(t)E^{\ast}}{S^{\ast}I^{\ast}E_{1}(t)}\right)  \right] \\
&  -\frac{S^{\ast}\beta bA^{\ast}}{S^{0}}\left[  g\left(  \frac{E_{1}%
(t)A^{\ast}}{E^{\ast}A_{1}(t)}\right)  +g\left(  \frac{S_{1}(t)A_{1}%
(t)E^{\ast}}{S^{\ast}A^{\ast}E_{1}(t)}\right)  \right] \\
\leq &  0.
\end{split}
\label{eq15}%
\end{equation}

From \eqref{eq14} and \eqref{eq15}, it follows $\omega(\phi)\subseteq
\mathfrak{D}$. Thereupon, if $\mathcal{R}_{c}>1$, $L$ is a Lyapunov function
on $\{v(t):t\geq1\}\subseteq\mathfrak{D}$. By \cite[Lemma 4.1]{Bai21}, we can
get that for any $\psi=(\psi_{1},\psi_{2},\psi_{3},\psi_{4},\psi_{5},\psi
_{6})^{T}\in\omega(\phi)$, there hold $\psi_{1}=S^{\ast}$, $\psi_{2}I^{\ast
}=E^{\ast}\psi_{3}$, $\psi_{2}A^{\ast}=E^{\ast}\psi_{4}$. Let $v(t)$ be the
solution of model \eqref{mod2} with any $\psi\in\omega(\phi).$ Accordingly,
the invariance of $\omega(\phi)$ implies that $S(t)=S^{\ast}$, $E^{\ast}%
I_{1}(t)=E_{1}(t)I^{\ast}$ and $E_{1}(t)A^{\ast}=E^{\ast}A_{1}(t)$ for all
$t\in\mathbb{R}$. We thus can obtain%
\[
E^{\ast}\dot{I}_{1}(t)=(pcE^{\ast}-B_{1}I^{\ast})E_{1}(t)=0\text{, }\forall
t\in\mathbb{R,}%
\]
thereby we have $\dot{I}_{1}(t)=\dot{E}_{1}(t)=\dot{A}_{1}(t)=0$ for
$t\in\mathbb{R.}$ Then the functions $E_{1}(t),$ $I_{1}(t)$ and $A_{1}(t)$ are
all constant functions on $\mathbb{R}$. Model \eqref{mod2} and the
invariance of $\omega(\phi)$ imply that both $Q(t)$ and $R(t)$ are constant
functions. Accordingly, $v(t)$ is a positive equilibrium of model
\eqref{mod2}. From \cite[Lemma 3.1]{Guo22}, it follows $v(t)=V^{\ast}$ for all
$t\in\mathbb{R}$. Consequently, $\omega(\phi)=\{V^{\ast}\}$ for $\mathcal{R}%
_{c}>1$, and hence $W^{s}(V^{\ast})=\Omega,$ where $W^{s}(V^{\ast})$ is the
stable set of $V^{\ast}$ for model \eqref{mod2}. By using the similar argument
as in Theorem 2.1, we can gain that $V^{\ast}$ is locally asymptotically stable
for model \eqref{mod2}. It follows easily from Theorem \ref{thm1} (or
\cite[Theorem 4.2]{Guo22}) that $\omega(\varphi)\cap W^{s}(V^{\ast}%
)\neq\emptyset$. Therefore, \cite[Theorem 1.2]{Thieme92} implies that
$\omega(\varphi)=\{V^{\ast}\}.$

\section{Dynamics of short-term COVID-19 model}

At present, the COVID-19 can be cleared out in a short term in some cities. Cui et
al. \cite{Cui20} proposed a short-term COVID-19 model with the contact rate
associated with real-time data on confirmed cases in Wuhan and Guangzhou at the year 2020. In their model, the
population birth and death rates were no longer considered. Therefore, we can
assume the total population $N(t)$ (defined by \eqref{nstt}) is a positive
constant, denoted it as $\mathcal{N}$. That is, for a short-term COVID-19 model, we can
take $\lambda=d=0$  in model \eqref{mod1} as follows
\begin{equation}%
\begin{split}
\dot{S}(t)  &  =-\beta\frac{S(t)}{N(t)}(aE(t)+I(t)+bA(t)),\\
\dot{E}(t)  &  =\beta\frac{S(t)}{N(t)}(aE(t)+I(t)+bA(t))-cE(t),\\
\dot{I}(t)  &  =pcE(t)-(q_{1}+r_{1})I(t),\\
\dot{A}(t)  &  =(1-p)cE(t)-(q_{2}+r_{2})A(t),\\
\dot{Q}(t)  &  =q_{1}I(t)+q_{2}A(t)-r_{3}Q(t),\\
\dot{R}(t)  &  =r_{1}I(t)+r_{2}A(t)+r_{3}Q(t).
\end{split}
\label{mod7}%
\end{equation}

In the following, we will study the dynamical behavior of this model.
Specifically, we calculate the control reproduction number, analyze the stability of the equilibria and obtain the expression of the final size
of the short-term model.

All solutions of the model with the nonnegative initial values exist, which are
unique, nonnegative and satisfy $N(t)=\mathcal{N}$. The model \eqref{mod7}
has multiple COVID-19-free equilibria $\tilde{V}=\left(  \tilde
{S},0,0,0,0,\mathcal{N}-\tilde{S}\right)  ^{T}$, $0\leq\tilde{S}%
\leq\mathcal{N}$, but there is no pandemic equilibrium. Using the method in
\cite{Driessch13}, the control reproduction number of \eqref{mod7} can be
taken by
\begin{equation}
\mathcal{R}_{c}=\frac{\tilde{S}}{\mathcal{N}}\left[  \frac{a\beta}{c}%
+\frac{p\beta}{q_{1}+r_{1}}+\frac{b\beta(1-p)}{q_{2}+r_{2}}\right]  .
\label{scrn}%
\end{equation}
Thanks to the control reproduction number $\mathcal{R}_{c}$, there are
three reasons for the spread of the pandemic, namely exposed individuals,
symptomatic infected individuals and asymptomatic infected individuals. In
the same line as the long-term model, the partial derivatives of
\eqref{scrn} with respect to $b$, $1-p$, $q_{1}$ and $q_{2}$ can be taken as follows
\[
\frac{\partial\mathcal{R}_{c}}{\partial b}=\frac{\tilde{S}}{\mathcal{N}}%
\frac{\beta(1-p)}{q_{2}+r_{2}},~\frac{\partial\mathcal{R}_{c}}{\partial
(1-p)}=\frac{\tilde{S}}{\mathcal{N}}\left(  \frac{b\beta}{q_{2}+r_{2}}%
-\frac{\beta}{q_{1}+r_{1}}\right)  ,~\frac{\partial\mathcal{R}_{c}}{\partial
q_{1}}=\frac{\tilde{S}}{\mathcal{N}}\frac{-p\beta}{(q_{1}+r_{1})^{2}}%
,~\frac{\partial\mathcal{R}_{c}}{\partial q_{2}}=\frac{\tilde{S}}{\mathcal{N}%
}\frac{-b\beta(1-p)}{(q_{2}+r_{2})^{2}}.
\]
The conclusions obtained here are consistent with the long-term model, which
will be explained in details in Sections \ref{622} and \ref{ScsSA}.

\subsection{Stability analysis}

It is not difficult to find that the set $\Gamma=\left\{  \phi\in
\mathbb{R}_{+}^{6}:\sum_{i=1}^{6}\phi_{i}=\mathcal{N}\right\}  $ is positively
invariant for model \eqref{mod7}, which is well-posed and dissipative in
$\Gamma.$ Next, we will discuss the dynamics of model \eqref{mod7} in $\Gamma
$. Let $U\left(  t\right)  :=\left(  S(t),E(t),I(t),A(t),Q(t),R(t)\right)
^{T}$ be the solution of model \eqref{mod7} through any $\phi\in\Gamma$, and
$\omega\left(  \phi\right)  $ be the $\omega$-limit set of $\phi$ for
$U\left(  t\right)  $. Then we have the following result.

\begin{theorem}
\label{thmS1} It holds that $\omega\left(  \phi\right)  \subseteq\left\{
\phi\in\mathbb{R}_{+}^{6}:\phi_{1}+\phi_{6}=\mathcal{N},\phi_{i}=0,i=2,3,4,5\right\}\subseteq \Gamma.$
\end{theorem}

\proof  From the first two equations of model \eqref{mod7}, it follows
\begin{equation}
\label{dot1}\dot{S}(t)+\dot{E}(t)=-cE(t).
\end{equation}
We have $\lim_{t\rightarrow\infty}(\dot{S}(t)+\dot{E}(t))=0$ by using the
similar method as in \cite{Lv20}, and hence $\lim_{t\rightarrow\infty}E(t)=0$.
Thus, it follows from model \eqref{mod7} that
\[
\lim_{t\rightarrow\infty}I(t)=\lim_{t\rightarrow\infty}A(t)=\lim
_{t\rightarrow\infty}Q(t)=0.
\]
Since $S(t)+E(t)+I(t)+A(t)+Q(t)+R(t)=\mathcal{N}$ for all $t\geq0,$ we obtain
$\lim_{t\rightarrow\infty}(S(t)+R(t))=\mathcal{N}$.

\begin{theorem}
\label{thmS2} The equilibrium $\hat{V}=\left(  0,0,0,0,0,\mathcal{N}\right)
^{T}$ is stable in $\Gamma$.
\end{theorem}

\proof Define the following function $V$ on $\Gamma,$
\[
V(\phi)=\phi_{1}+\phi_{2}+\phi_{3}+\phi_{4}+\phi_{5}+\mathcal{N}-\phi_{6}.
\]
It easily follows that $V$ is a positive definite function with respect to
$\hat{V}$. Consequently, the derivative of $V$ along the solution $U(t)$ is%
\[
\dot{V}(U(t)) =\dot{S}(t)+\dot{E}(t)+\dot{I}(t)+\dot{A}(t)+\dot{Q} (t)-\dot
{R}(t)=-2\left(  r_{1}I(t)+r_{2}A(t)+r_{3}Q(t)\right)  \leq0.
\]
Thus $\hat{V}$ is stable.

\begin{theorem}
\label{thmS3} The equilibrium $\tilde{V}=\left(  \tilde{S},0,0,0,0,\mathcal{N}%
-\tilde{S}\right)  $ $(\tilde{S}\in(0,\mathcal{N}])$ is stable if
$\mathcal{R}_{c}\leq1$ and unstable if $\mathcal{R}_{c}>1$ in $\Gamma
_{1}=\left\{  \phi\in\Gamma:\phi_{1}>0\right\}  $.
\end{theorem}

\proof Obviously, $\Gamma_{1}$ is a positive invariant set of
model \eqref{mod7}. Let $\Gamma_{2}=\left\{  \chi=(\chi_{1},\chi_{2},\chi
_{3},\chi_{4})^{T}\in\mathbb{R}_{+}^{4}:\chi_{1}>0\right\}  .$ Then
$\Gamma_{2}$ is positively invariant for the following model
\begin{equation}%
\begin{split}
\dot{S}(t) &  =-\beta\frac{S(t)}{\mathcal{N}}(aE(t)+I(t)+bA(t)),\\
\dot{E}(t) &  =\beta\frac{S(t)}{\mathcal{N}}(aE(t)+I(t)+bA(t))-cE(t),\\
\dot{I}(t) &  =pcE(t)-(q_{1}+r_{1})I(t),\\
\dot{A}(t) &  =(1-p)cE(t)-(q_{2}+r_{2})A(t),
\end{split}
\label{mods}%
\end{equation}
Clearly, $X^{0}=\left(  \tilde{S},0,0,0\right)  $ is an equilibrium of model
\eqref{mods}. The function $V$ is defined as follows
\[
V(\chi)=\mathcal{N}\left(  \frac{\chi_{1}}{\tilde{S}}-1-\ln\frac{\chi_{1}%
}{\tilde{S}}\right)  +\frac{\mathcal{N}}{\tilde{S}}\chi_{2}+\frac{\beta}%
{q_{1}+r_{1}}\chi_{3}+\frac{\beta b}{q_{2}+r_{2}}\chi_{4},\text{ }\chi
\in\Gamma_{2}.\label{V51}%
\]
Observe that $V$ is a positive definite function with respect to $X^{0}$. Let
$u\left(  t\right)  :=\left(  S(t),E(t),I(t),A(t)\right)  ^{T}$ be the
solution of model \eqref{mods} through any $\chi\in\Gamma_{2}$. Then for
$\mathcal{R}_{c}\leq1,$ the derivative of $V$ along $u\left(
t\right)  $ is
\[%
\begin{split}
\dot{V}\left(  u\left(  t\right)  \right)   &  =\mathcal{N}\left(  \frac
{1}{\tilde{S}}-\frac{1}{S(t)}\right)  \dot{S}(t)+\frac{\mathcal{N}}{\tilde{S}%
}\dot{E}(t)+\frac{\beta}{q_{1}+r_{1}}\dot{I}(t)+\frac{\beta b}{q_{2}+r_{2}%
}\dot{A}(t)\\
&  =\beta\left(  aE(t)+I(t)+bA(t)\right)  -\frac{\mathcal{N}}{\tilde{S}%
}cE(t)+\frac{\beta}{q_{1}+r_{1}}pcE(t)-\beta I(t)+\frac{\beta b}{q_{2}+r_{2}%
}(1-p)cE(t)-\beta bA(t)\\
&  =cE(t)\frac{\mathcal{N}}{\tilde{S}}\left(  R_{c}-1\right)  \leq0.
\end{split}
\label{V52}%
\]
In consequence, $X^{0}$ is stable in $\Gamma_{2}$. By using the method in
\cite{Bai21}, we can obtain $\tilde{V}$ is stable in $\Gamma_{1}$.

Next, we show that $\tilde{V}$ is unstable for $\mathcal{R}_{c}>1$. In fact,
the characteristic equation of the linearized system corresponding to model
\eqref{mod7} at $\tilde{V}$ is given by
\[
f(\Lambda)=\Lambda^{2}\left(  \Lambda+r_{3}\right)  \left(  \Lambda^{3}%
+a_{1}\Lambda^{2}+a_{2}\Lambda+a_{3}\right)  ,
\]
where
\begin{align*}
a_{1}  &  =c(1-\mathcal{R}_{c})+(q_{1}+r_{1})+(q_{2}+r_{2})+\frac{\tilde{S}%
}{\mathcal{N}}\frac{pc\beta}{q_{1}+r_{1}}+\frac{\tilde{S}}{\mathcal{N}}%
\frac{bc\beta(1-p)}{q_{2}+r_{2}},\\
a_{2}  &  =(q_{1}+r_{1}+q_{2}+r_{2})c(1-\mathcal{R}_{c})+\frac{\tilde{S}%
}{\mathcal{N}}\frac{(q_{2}+r_{2})pc\beta}{(q_{1}+r_{1})}+\frac{\tilde{S}%
}{\mathcal{N}}\frac{(q_{1}+r_{1})bc\beta(1-p)}{(q_{2}+r_{2})}+(q_{1}%
+r_{1})(q_{2}+r_{2}),\\
a_{3}  &  =(q_{1}+r_{1})(q_{2}+r_{2})c(1-\mathcal{R}_{c})<0.
\end{align*}
Therefore, $f(\Lambda)=0$ has a positive root. As a result, $\tilde{V}$ is unstable.

\subsection{The final COVID-19 size}

Based on model \eqref{mod7}, we can get the final COVID-19 size,
which represents the percentage of the infected population. The
initial values of model \eqref{mod7} are $S(0)=\mathcal{N}%
-E(0)-I(0)-A(0)-Q(0)>0$, $E(0)\geq0$, $I(0)\geq0$, $A(0)\geq0$, $Q(0)\geq0$,
$R(0)=0$. From Theorem \ref{thmS1}, it follows that
\[
\lim_{t\rightarrow\infty}E(t)=\lim_{t\rightarrow\infty}I(t)=\lim
_{t\rightarrow\infty}A(t)=\lim_{t\rightarrow\infty}Q(t)=0. \label{inf}%
\]
The existence of $S(\infty):=\lim_{t\rightarrow\infty}S(t)$ can be
derived from the boundedness and monotonicity of $S(t)$. By \eqref{dot1},
there holds
\[
\int_{0}^{\infty}E(t)dt=\frac{E(0)+S(0)-S(\infty)}{c}. \label{dotE}%
\]
By the first three equations of model \eqref{mod7}, we can obtain
\[
\int_{0}^{\infty}I(t)dt=\frac{pE(0)+pS(0)+I(0)-pS(\infty)}{q_{1}+r_{1}}.
\label{dotI}%
\]
From the first four equations of model \eqref{mod7}, we can derive
\[
\int_{0}^{\infty}A(t)dt=\frac{(1-p)S(0)+(1-p)E(0)+A(0)-(1-p)S(\infty)}%
{q_{2}+r_{2}}. \label{dotA}%
\]
It follows from the first equation of model \eqref{mod7} that
\[
S(\infty)=S(0)e^{-\frac{\beta}{\mathcal{N}}\int_{0}^{\infty}\left(
aE(t)+I(t)+bA(t)\right)  dt}>0
\]
and
\[%
\begin{split}
\ln\frac{S(0)}{S(\infty)}={}  &  \frac{\beta a}{\mathcal{N}}\int_{0}^{\infty
}E(t)dt+\frac{\beta}{\mathcal{N}}\int_{0}^{\infty}I(t)dt+\frac{\beta
b}{\mathcal{N}}\int_{0}^{\infty}A(t)dt\\
={}  &  \frac{\beta a}{\mathcal{N}}\frac{S(0)+E(0)-S(\infty)}{c}+\frac{\beta
}{\mathcal{N}}\frac{pS(0)+pE(0)+I(0)-pS(\infty)}{q_{1}+r_{1}}\\
&  +\frac{\beta b}{\mathcal{N}}\frac{(1-p)S(0)+(1-p)E(0)+A(0)-(1-p)S(\infty
)}{q_{2}+r_{2}}\\
={}  &  \frac{S(0)+E(0)-S(\infty)}{\tilde{S}}R_{c}+\frac{\beta I(0)}%
{\mathcal{N}(q_{1}+r_{1})}+\frac{\beta bA(0)}{\mathcal{N}(q_{2}+r_{2})},
\end{split}
\]
where $\tilde{S}>0$. Thus, the total infected population is
\[%
\begin{split}
Z  &  =S(0)-S(\infty)\\
&  =S(0)\left\{  1-e^{-\left[\frac{S(0)+E(0)-S(\infty)}{\tilde{S}}%
R_{c}+\frac{\beta I(0)}{\mathcal{N}(q_{1}+r_{1})}+\frac{\beta bA(0)}%
{\mathcal{N}(q_{2}+r_{2})}\right]  }\right\} \\
&  =S(0)\left\{  1-e^{-\left[  \frac{Z+E(0)}{\tilde{S}}R_{c}+\frac{\beta
I(0)}{\mathcal{N}(q_{1}+r_{1})}+\frac{\beta bA(0)}{\mathcal{N}(q_{2}+r_{2}%
)}\right]  }\right\}  .
\end{split}
\]
Therefore, the expression of the final COVID-19 size is $\mathcal{F}:=\frac{Z}{\mathcal{N}}$.

\section{Case study}

To illustrate the impact of quarantine measures and asymptomatic transmission
on COVID-19, we will use long-term model and short-term model to demonstrate
COVID-19 transmission in India and Nanjing, respectively. Besides,
sensitivity analysis of $\mathcal{R}_{c}$ will be conducted in these two
areas, which provides theoretical basis for proposing COVID-19 prevention and
control measures. It is worth mentioning that in the numerical analysis of
short-term model, we will choose $\tilde{S}=\mathcal{N}$ for most common
cases.

\subsection{Case study in long-term COVID-19 model}

\begin{figure}[ptb]
\centering
\begin{minipage}[t]{0.45\linewidth}
\centering
\includegraphics[scale=0.52]{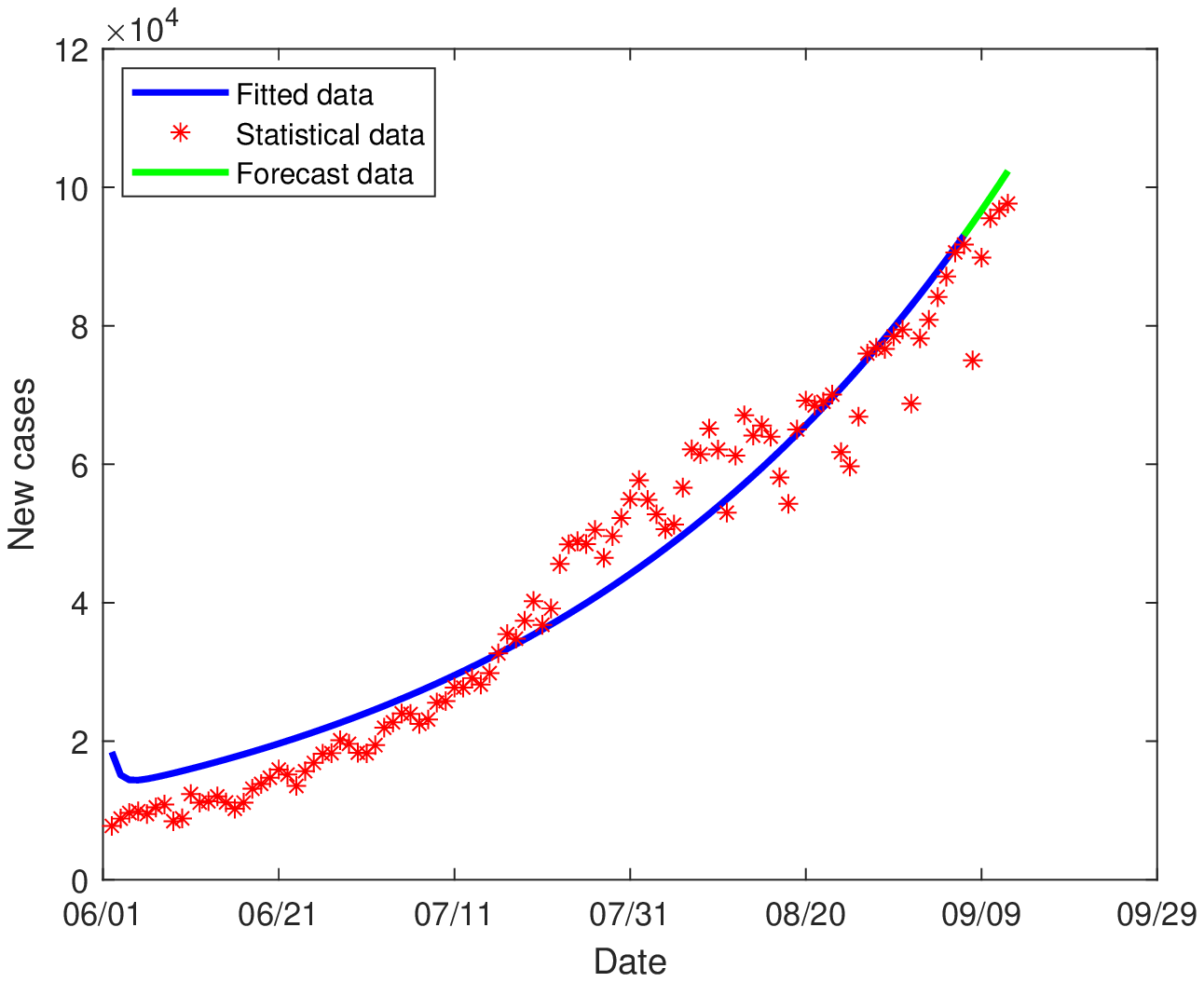}
\caption{The daily new cases in India.}\label{fig6}
\end{minipage}\begin{minipage}[t]{0.45\linewidth}
\centering
\includegraphics[scale=0.52]{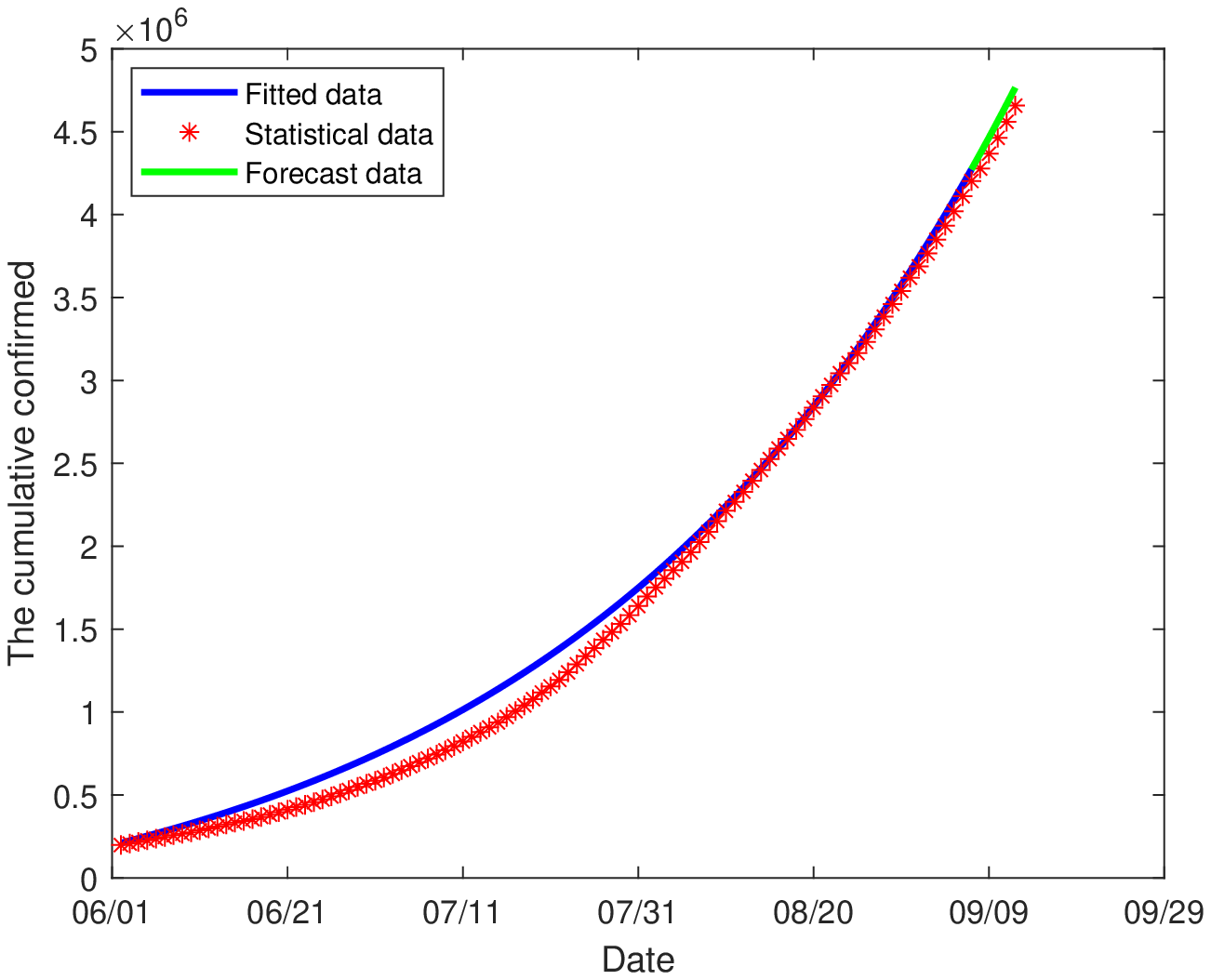}
\caption{Cumulative confirmed cases in India.}\label{fig7}
\end{minipage}
\end{figure}

\begin{figure}[ptb]
\centering
\begin{minipage}[t]{0.45\linewidth}
\includegraphics[scale=0.52]{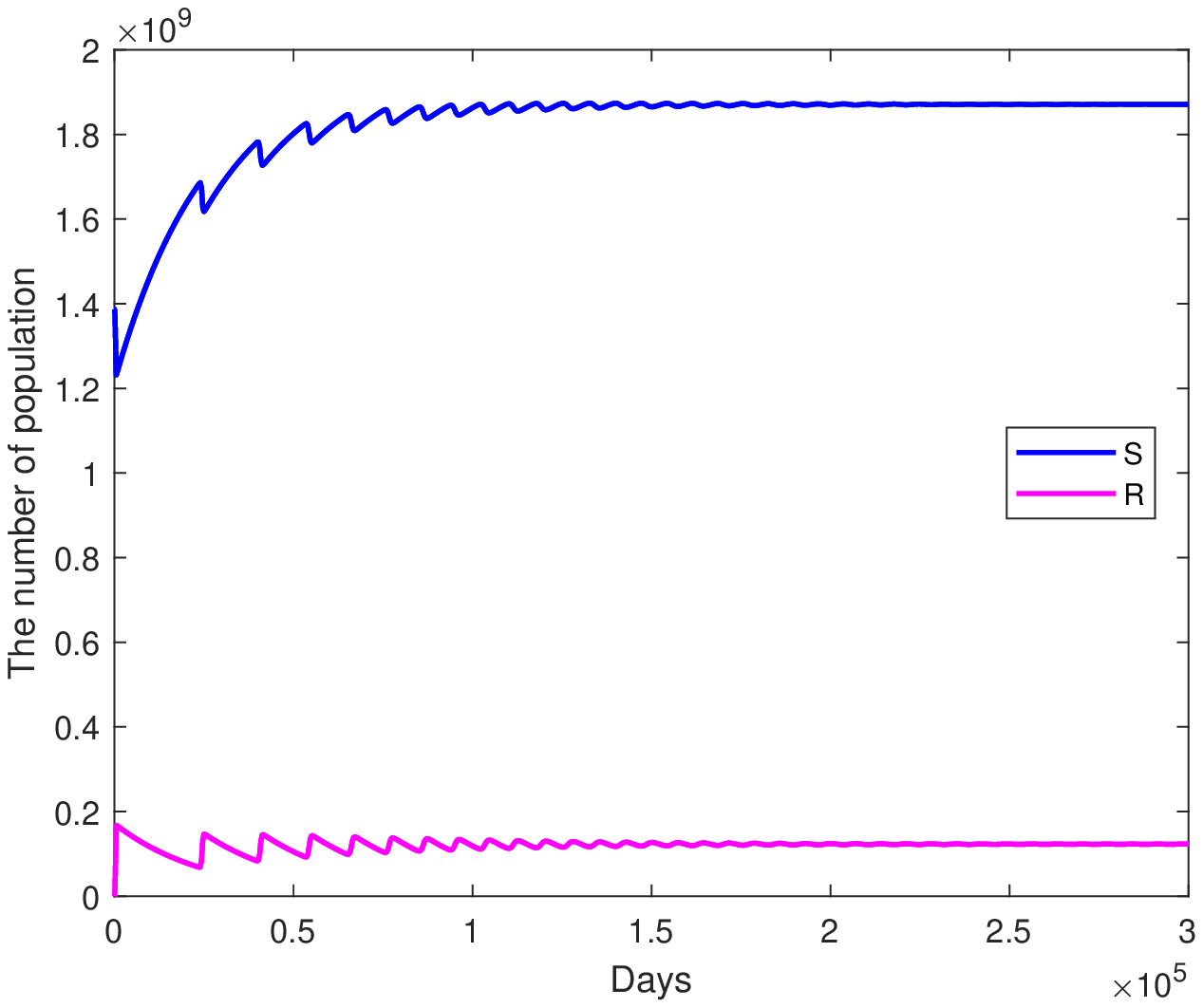}
\caption{Time evolutions of population in \\ $S$ and $R$ compartments in India.}\label{fig8}
\end{minipage}\begin{minipage}[t]{0.45\linewidth}
\centering
\includegraphics[scale=0.52]{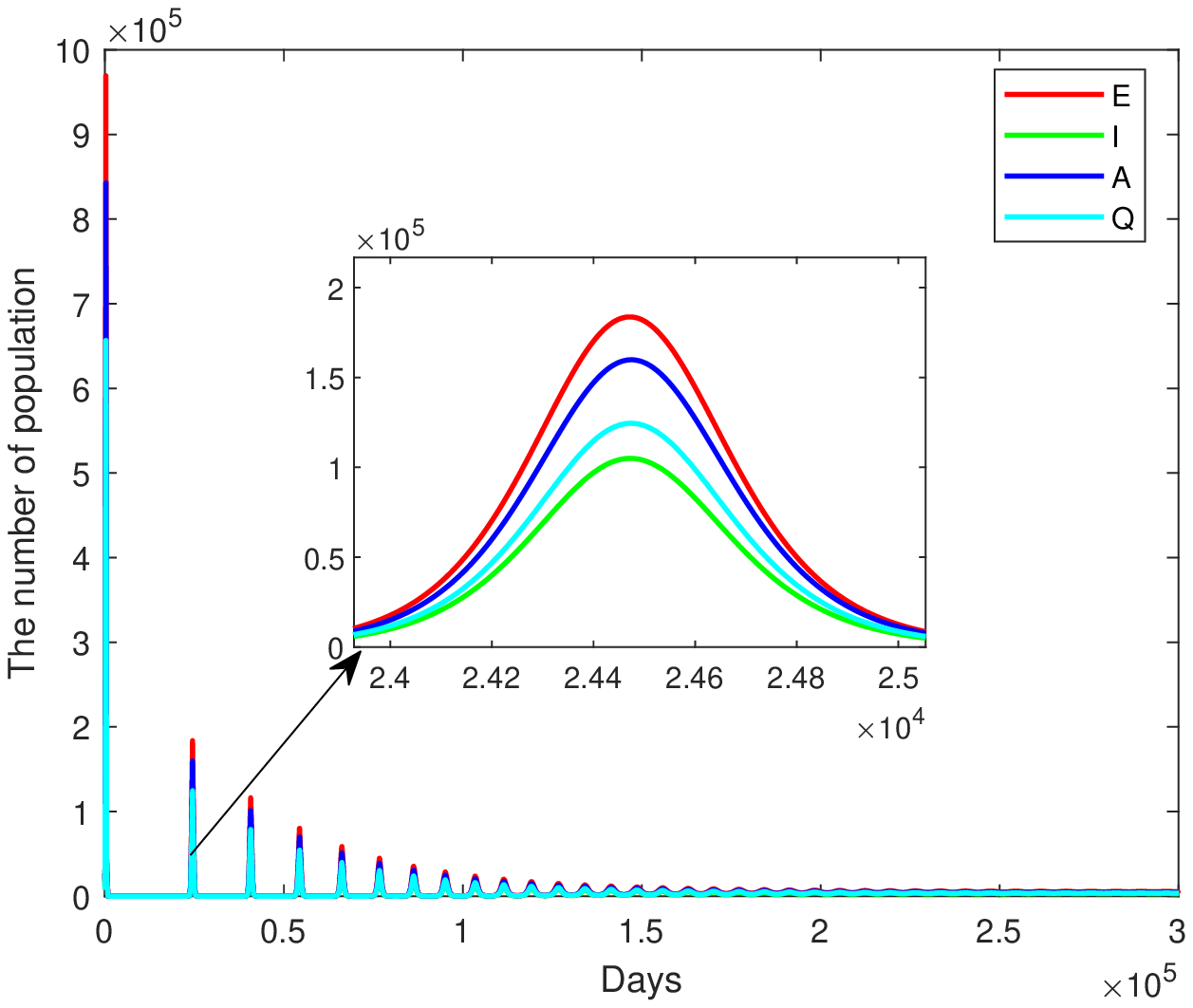}
\caption{Time evolutions of population in \\ different compartments in India.}\label{fig9}
\end{minipage}
\end{figure}

COVID-19 cases emerged on February 2020 in India and as of now new cases is
still reported every day \cite{Indiadata}. At this point, it is more
reasonable to use the long-term model \eqref{mod1} to illustrate the
transmission dynamics of COVID-19 in India. In \cite{Yuan21}, we know that the
population of India is about $N=1386750000$ on May 2020. The number of daily
new cases and cumulative number of confirmed cases in India from June 1, 2020
to September 11, 2020 are obtained from \cite{Indiadata}. The number of daily
new cases from June 1, 2020 to September 6, 2020 have been used to estimate
the parameters in Tab. \ref{tab:parameter}.
 \begin{table}[ptbh]
\caption{Parameter estimates of COVID-19 model \eqref{mod1}.}%
\label{tab:parameter}%
\centering
\begin{tabular}
[c]{ccc}%
\toprule[0.6mm] Parameter & Parameter value of
India & Source\\\hline
$\lambda$ & $7.7575\times10^{4}$ &  \cite{Yuan21}\\
$d $ & $3.8905\times10^{-5}$ &  \cite{Yuan21}\\
$\beta$ & 0.3717 & Estimated\\
$a $  & 0.4898 & Estimated\\
$b $ & 1.7281 & Estimated\\
$c $ & 0.8945 & Estimated\\
$p $ & 0.6937 & Estimated\\
$q_{1} $  & 0.8296 & Estimated\\
$q_{2} $ & 0.1947 & Estimated\\
$r_{1} $ & 0.2565 & Estimated\\
$r_{2} $ & 0.1201 & Estimated\\
$r_{3} $ & 0.9495 & Estimated\\\hline\hline
$S(0) $ & $1.3867\times10^{9}$ & Estimated\\
$E(0) $ & $1.9727\times10^{4}$ & Estimated\\
$I(0) $ & $1.8179\times10^{4}$ & Estimated\\
$A(0) $ & $1.7174\times10^{4}$ & Estimated\\
$Q(0) $ & $7.5236\times10^{3}$ & Estimated\\
$R(0) $ & $1.3321\times10^{4}$ & Estimated\\
\bottomrule[0.6mm] &  &
\end{tabular}
\end{table}
Using the parameter values of
India in Tab. \ref{tab:parameter}, we have fitted the daily new cases and
cumulative confirmed cases (blue curve) and compared them with the statistical
data (red curve), as shown in Figs. \ref{fig6} and \ref{fig7}. Moreover, we
predicted daily new cases and cumulative confirmed cases for the five days
from September 7, 2020 to September 11, 2020 (green curve), and compared them
with the corresponding statistical values. The control reproduction number $\mathcal{R}_{c}$ in India have
been estimated to be $1.0657>1$. If the authorities do nothing to change the
existing situation, it can be seen from Figs. \ref{fig8} and \ref{fig9} that
the populations of all compartments in India will eventually tend to the
COVID-19 equilibrium $V^{\ast}$. That is, COVID-19 will persist, which is consistent with our theoretical result in
Theorem \ref{thm52}. 

\subsubsection{The impact of quarantine measures on COVID-19}

\begin{figure}[ptb]
\centering
\includegraphics[scale=0.52]{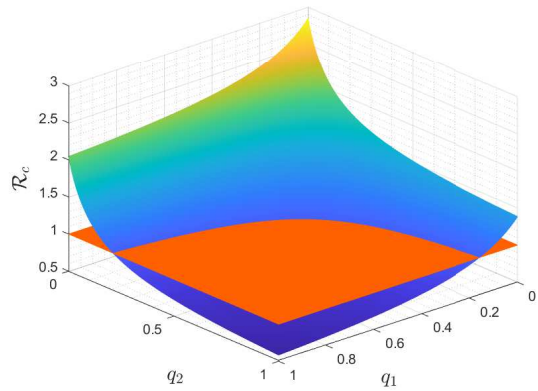} \caption{The relationship
among $\mathcal{R}_{c}$ and parameters $q_{1}$, $q_{2}$ in India.}%
\label{fig10}%
\end{figure}

\begin{figure}[ptb]
\centering
\begin{minipage}[t]{0.45\linewidth}
\centering
\includegraphics[scale=0.52]{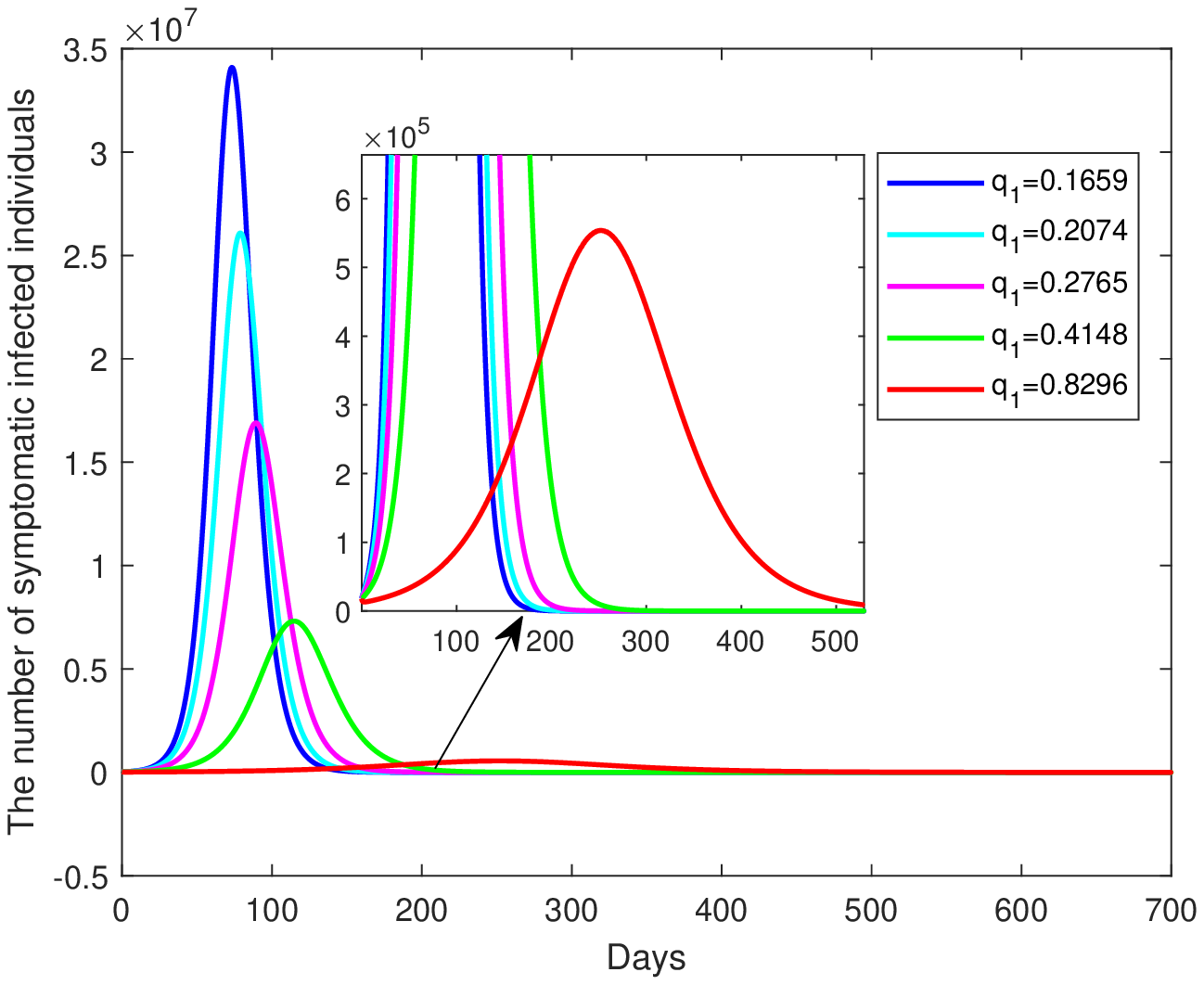}
\caption{The relationship between $I(t)$ and \\$q_{1}$ in model \eqref{mod1}.}\label{figL1}
\end{minipage}\begin{minipage}[t]{0.45\linewidth}
\centering
\includegraphics[scale=0.52]{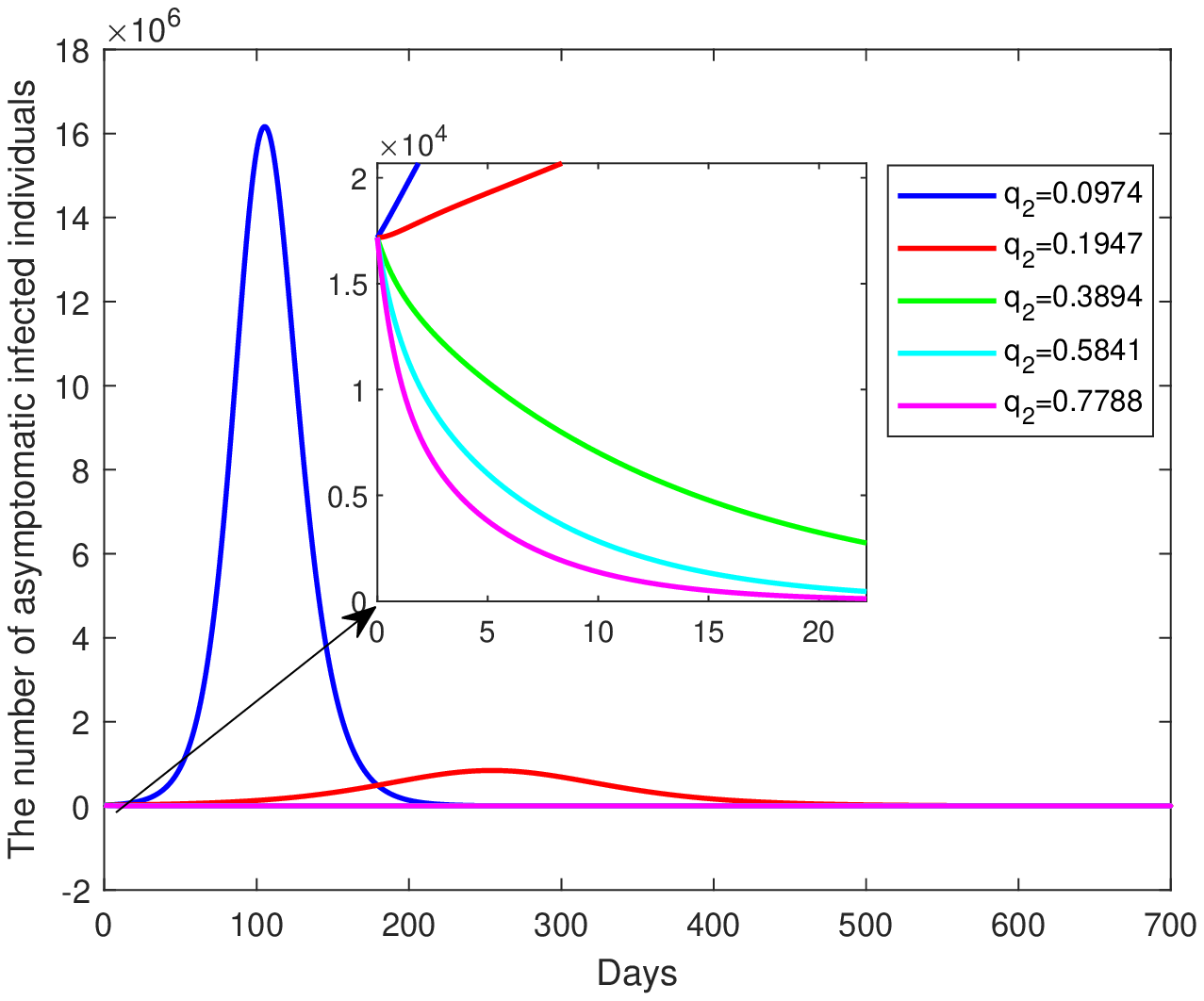}
\caption{The relationship between $A(t)$ and \\$q_{2}$ in model \eqref{mod1}.}\label{figL2}
\end{minipage}
\end{figure}

\begin{figure}[ptb]
\centering
\begin{minipage}[t]{0.45\linewidth}
\centering
\includegraphics[scale=0.52]{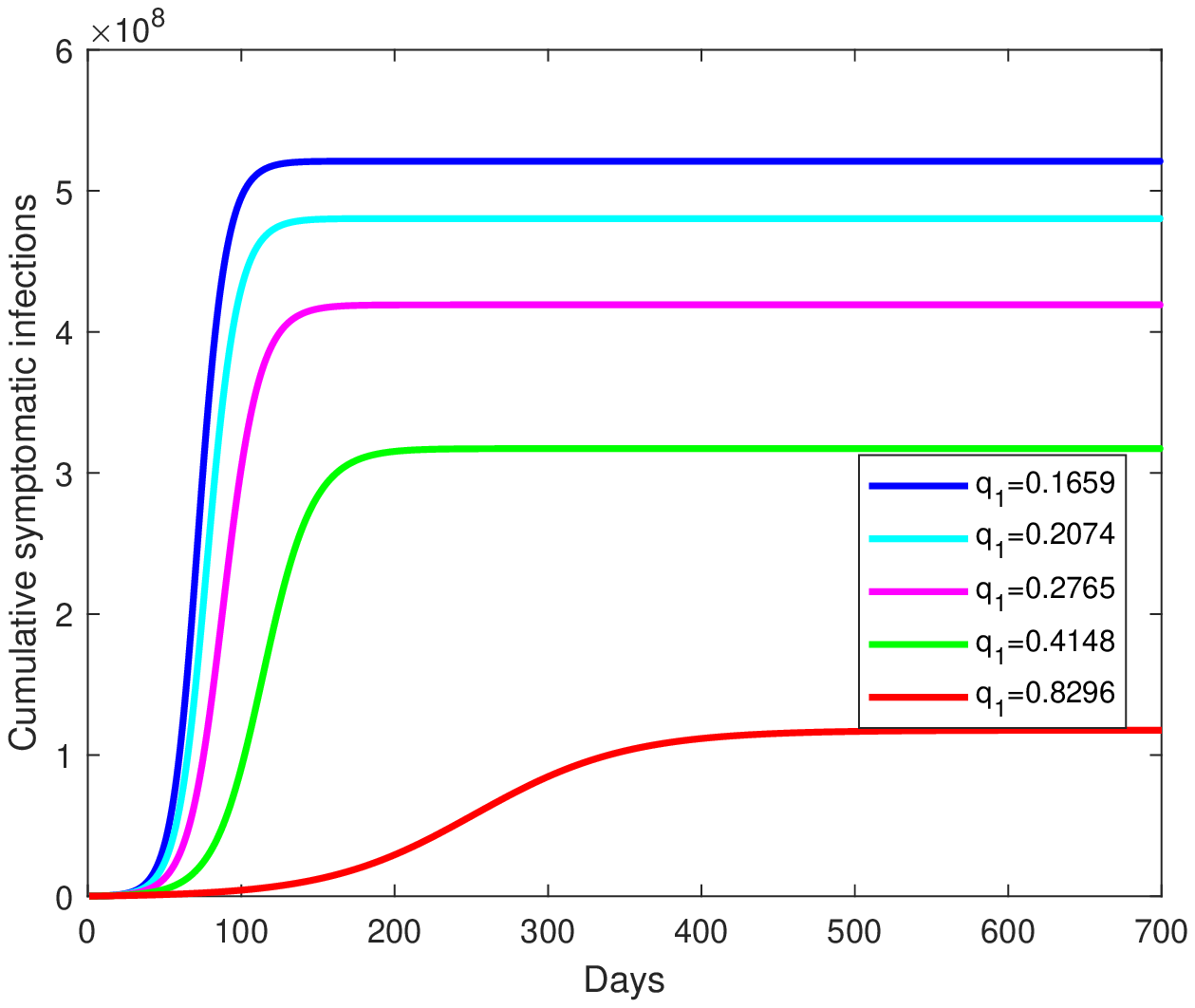}
\caption{The cumulative symptomatic \\infections in model \eqref{mod1}.}\label{figLC1}
\end{minipage}\begin{minipage}[t]{0.45\linewidth}
\centering
\includegraphics[scale=0.52]{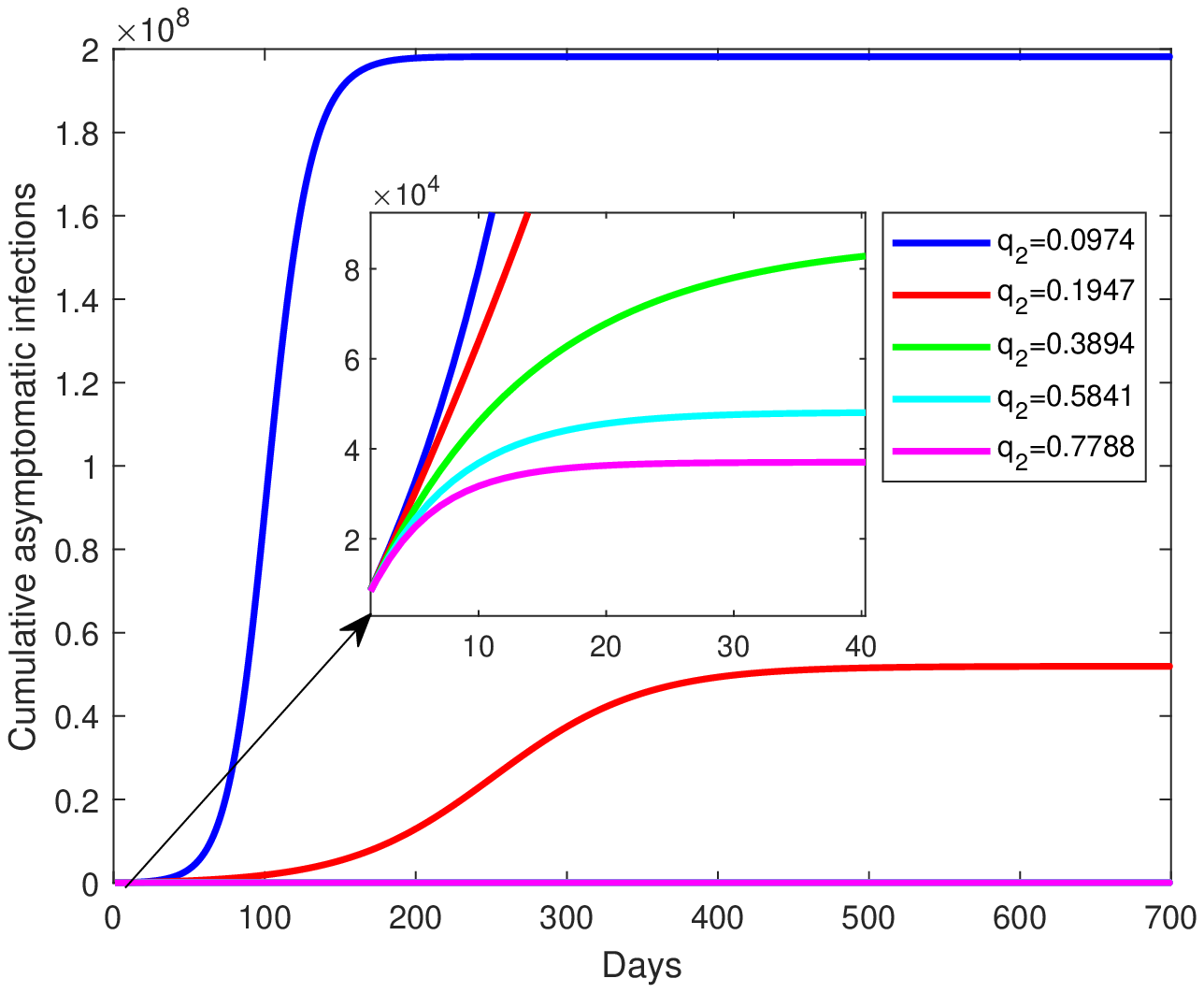}
\caption{The cumulative asymptomatic \\infections in model \eqref{mod1}.}\label{figLC2}
\end{minipage}
\end{figure}

It can be found from Fig. \ref{fig10} that as the values of $q_{1}$ and
$q_{2}$ increases, the value of $\mathcal{R}_{c}$ will decrease and can be
reduced to $\mathcal{R}_{c}<1$, which implies that COVID-19 will disappear,
see Theorem \ref{thm51}. Therefore, the strengthened quarantine measures can
be effective to control the COVID-19 pandemic. In Figs. \ref{figL1} and \ref{figLC1}, all parameter values except
$q_{1}$ remain unchanged. With the increase of value of $q_{1}$, both the peak
value of $I(t)$ and cumulative number of symptomatic COVID-19 infections will
decrease. In Figs. \ref{figL2} and \ref{figLC2}, only the value of parameter
$q_{2}$ is changed. Both the peak value of $A(t)$ and cumulative number of
asymptomatic COVID-19 infections decrease as the $q_{2}$ value increases.
This indicates that the peak value of infected individuals and the cumulative
confirmed cases will been reduced by strengthening quarantine measures.

\begin{figure}[ptb]
\centering
\includegraphics[scale=0.52]{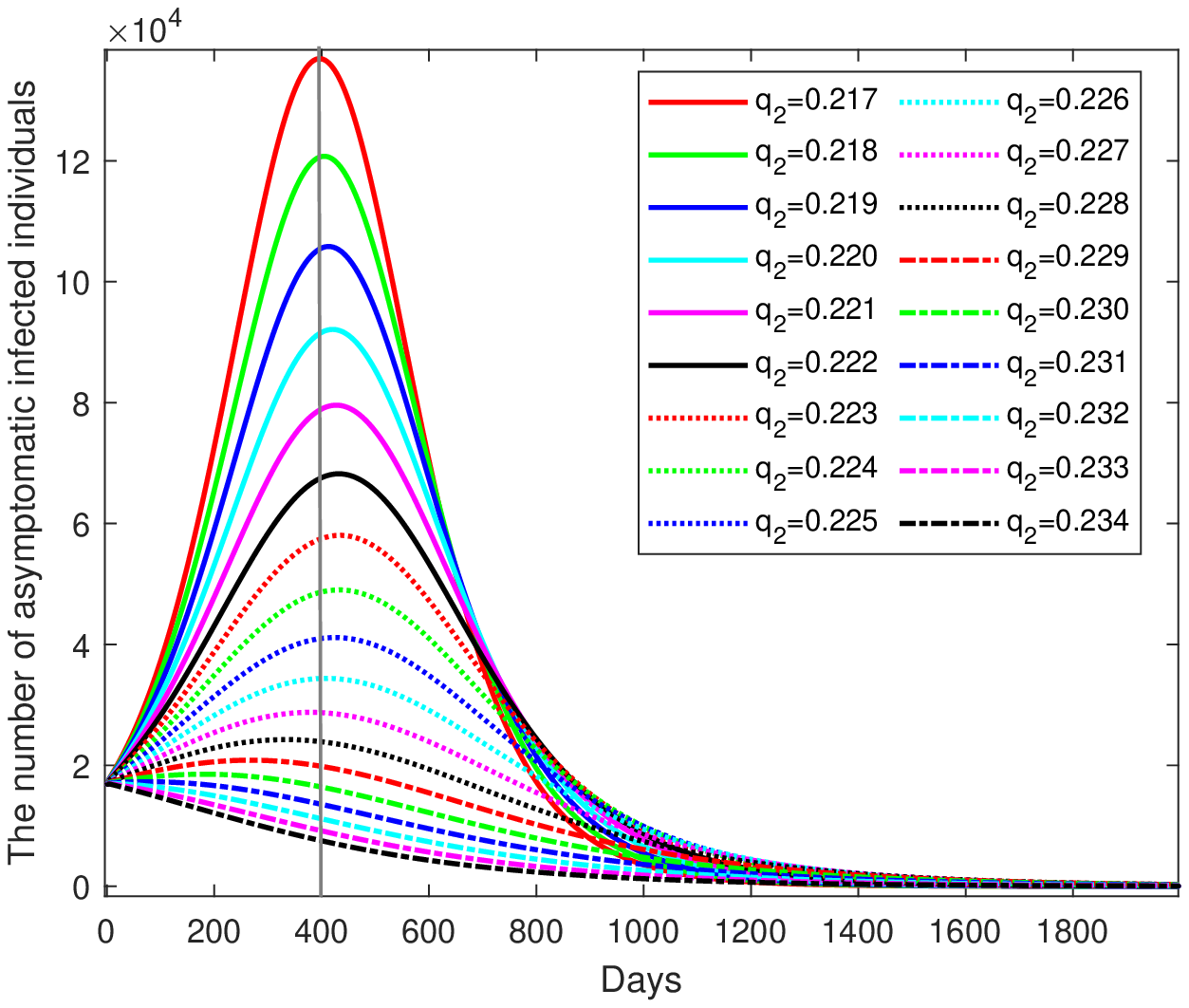}
\caption{The relationship between $A(t)$ and \\$q_{2}$ in model \eqref{mod1}.}\label{figcritical}
\end{figure}

In Fig. \ref{figL2}, when $q_{2}\leq0.1947$, days to reach the peak value of asymptomatic infected individuals gradually increase; when $q_{2}>0.1947$, this time gradually decreases. In order to show that this phenomenon is not a coincidence, we only change the quarantine rate of asymptomatic infected individuals and simulate the impact of different quarantine rates of asymptomatic infections on COVID-19 in India with other values unchanged. As shown in Fig. \ref{figcritical}, the quarantine rate of asymptomatic infections ranges from 0.217 to 0.234 in step of 0.001. Obviously, the smaller the step size is, the more accurate the critical value of $q_{2}$ will be. Tab. \ref{criticalV} more clearly shows the critical value $q_{2}=0.223$. It is easy to see that with the strengthening of quarantine measures, when $q_{2}\leq0.223$, the time to reach the peak value of asymptomatic infected individuals gradually increases; when $q_{2}>0.223$, this time gradually decreases.


\begin{table}[ptb]
\caption{Relationship between $q_{2}$ and days to reach the peak value of asymptomatic infected individuals.}%
\label{criticalV}%
\centering
\begin{tabular}
[c]{cc}%
\toprule The value of parameter $q_{2}$ & Days to reach the peak value\\
\midrule
0.219 & 413\\
0.220 & 422\\
0.221 & 427\\
0.222 & 433\\
0.223 & 436\\
0.224 & 433\\
0.225 & 424\\
0.226 & 408\\
\bottomrule &
\end{tabular}
\end{table}

\subsubsection{The impact of asymptomatic infections on COVID-19}

\label{612} 
\begin{figure}[ptb]
\centering
\begin{minipage}[t]{0.45\linewidth}
\centering
\includegraphics[scale=0.52]{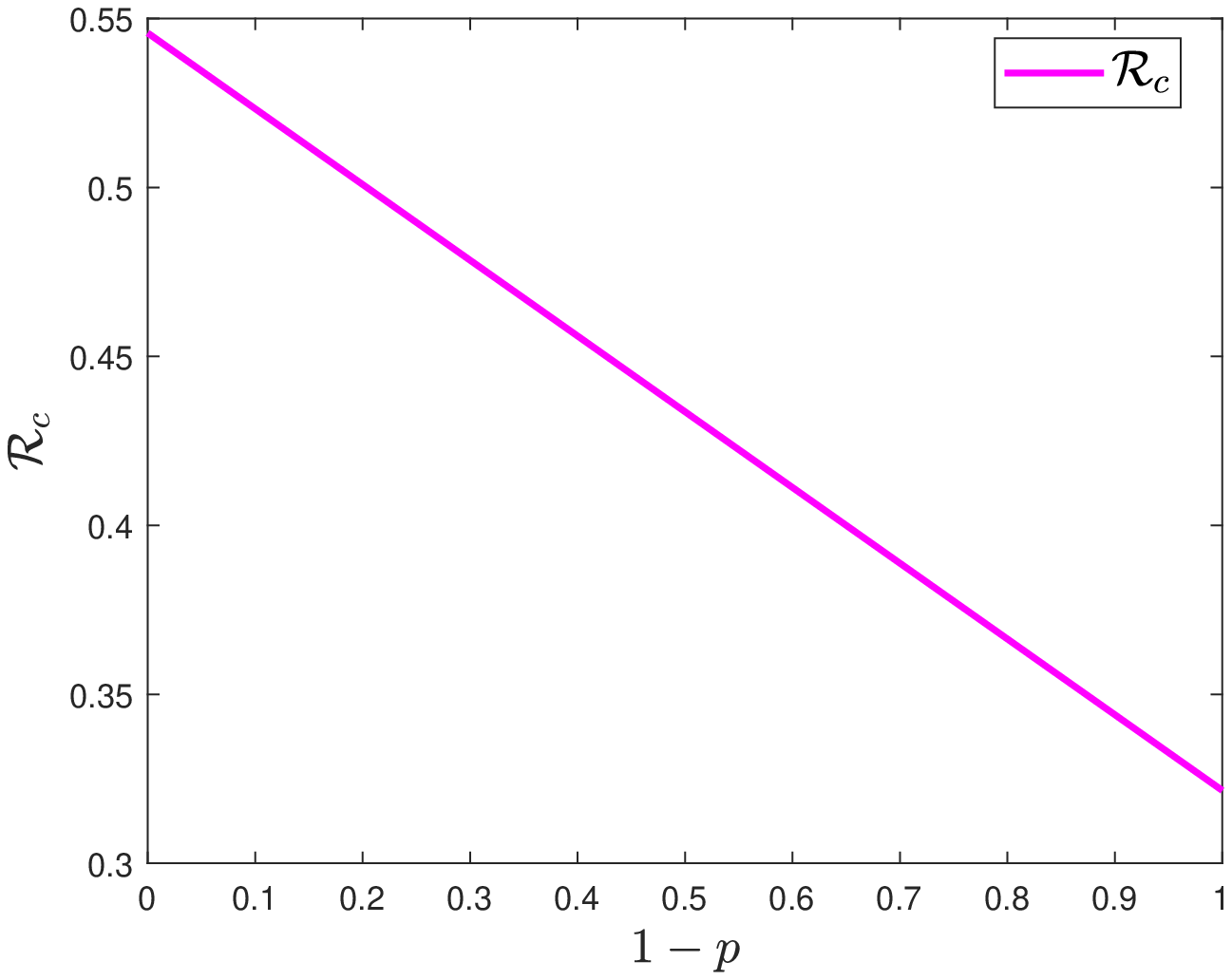}
\caption{The relationship between $\mathcal{R}_{c}$ and $1-p$.}\label{figLA1}
\end{minipage}\begin{minipage}[t]{0.45\linewidth}
\centering
\includegraphics[scale=0.52]{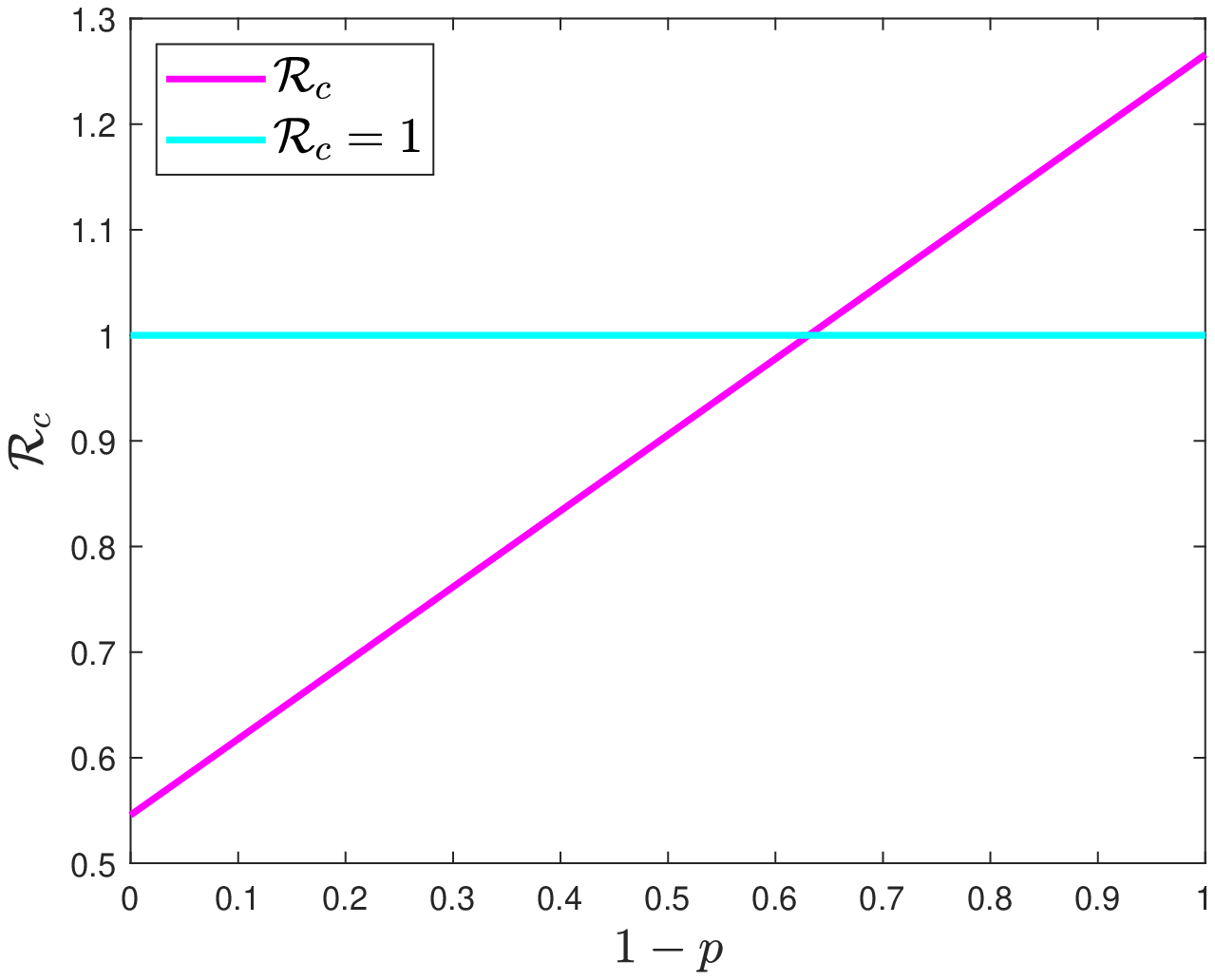}
\caption{The relationship between $\mathcal{R}_{c}$ and $1-p$.}\label{figLA2}
\end{minipage}
\end{figure}

\begin{figure}[ptb]
\centering
\includegraphics[scale=0.52]{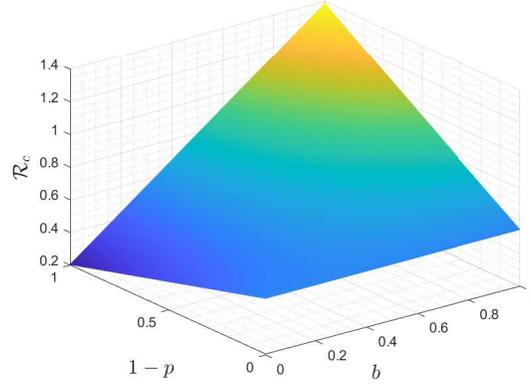} \caption{The relationship
among $\mathcal{R}_{c}$, $b$ and $1-p$ in India.}%
\label{figLA3}%
\end{figure}

Now, we analyze the impact of asymptomatic transmission on
COVID-19. Using the data of India in Tab. \ref{tab:parameter}, we can observe that if $b\leq
0.28991$, then $\mathcal{R}_{c}$ is negatively correlated with $1-p$ as shown in
Fig. \ref{figLA1}, while in case of $b> 0.28991$, $\mathcal{R}_{c}$ is positively
correlated with $1-p$ as shown in Fig. \ref{figLA2}. From Fig. \ref{figLA3},
it can be seen intuitively that $\mathcal{R}_{c}$ can change from greater than
1 to less than 1 when both $1-p$ and $b$ are changed. This suggests that
asymptomatic infections play a key role in the spread of COVID-19.


\subsubsection{Sensitivity analysis}

\label{LcsSA}The sensitivity index \cite{Chitnis08} of $\mathcal{R}_{c}$ with
respect to parameter $k$ is expressed as
\begin{equation}
\xi_{k}^{\mathcal{R}_{c}}=\frac{\partial\mathcal{R}_{c}}{\partial k}\cdot
\frac{k}{\mathcal{R}_{c}}. \label{tsi}%
\end{equation}
It can be seen from \eqref{crn} that $\mathcal{R}_{c}$ is affected by the
parameters $a$, $\beta$, $p$, $c$, $b$, $d$, $q_{1}$, $q_{2}$, $r_{1}$,
$r_{2}$. The values of parameters $p$, $c$ and $d$ are difficult to be changed
by artificial measures, so that they are fixed in the spread of COVID-19.
We will analyze the remaining seven parameters. According to \eqref{crn} and
\eqref{tsi}, we can obtain
\begin{align*}
\xi_{a}^{\mathcal{R}_{c}}  &  =\frac{\beta}{c+d}\cdot\frac{a}{\mathcal{R}_{c}%
}=\frac{aB_{1}B_{2}}{aB_{1}B_{2}+pcB_{2}+bc(1-p)B_{1}},\\
\xi_{\beta}^{\mathcal{R}_{c}}  &  =\left[  \frac{a}{c+d}+\frac{pc}{(c+d)B_{1}}
+\frac{bc(1-p)}{(c+d)B_{2}}\right]  \cdot\frac{\beta}{\mathcal{R}_{c}}=1,\\
\xi_{b}^{\mathcal{R}_{c}}  &  =\frac{\beta c(1-p)}{(c+d)B_{2}}\cdot\frac
{b}{\mathcal{R}_{c}}=\frac{cb(1-p)B_{1}} {aB_{1}B_{2}+pcB_{2}+cb(1-p)B_{1}},\\
\xi_{q_{1}}^{\mathcal{R}_{c}}  &  =-\frac{pc\beta}{(c+d)B_{1}^{2}}\cdot
\frac{q_{1}}{\mathcal{R}_{c}}=-\frac{pcq_{1}B_{2}}{B_{1}\left[  aB_{1}%
B_{2}+pcB_{2}+cb(1-p)B_{1}\right]  },\\
\xi_{q_{2}}^{\mathcal{R}_{c}}  &  =-\frac{bc\beta(1-p)}{(c+d)B_{2}^{2}}%
\cdot\frac{q_{2}}{\mathcal{R}_{c}}=-\frac{bc(1-p)q_{2}B_{1}}{B_{2}\left[
aB_{1}B_{2}+pcB_{2}+cb(1-p)B_{1}\right]  },\\
\xi_{r_{1}}^{\mathcal{R}_{c}}  &  =-\frac{pc\beta}{(c+d)B_{1}^{2}}\cdot
\frac{r_{1}}{\mathcal{R}_{c}}=-\frac{pcr_{1}B_{2}}{B_{1}\left[  aB_{1}%
B_{2}+pcB_{2}+cb(1-p)B_{1}\right]  },\\
\xi_{r_{2}}^{\mathcal{R}_{c}}  &  =-\frac{bc\beta(1-p)}{(c+d)B_{2}^{2}}%
\cdot\frac{r_{2}}{\mathcal{R}_{c}}=-\frac{bc(1-p)r_{2}B_{1}}{B_{2}\left[
aB_{1}B_{2}+pcB_{2}+cb(1-p)B_{1}\right]  }.\\
\end{align*}

By substituting the data of India in Tab. \ref{tab:parameter}, the sensitivity
index can be obtained as shown in Tab. \ref{tab3}. \begin{table}[ptb]
\caption{Sensitivity index of $\mathcal{R}_{c}$ with respect to parameters.}%
\label{tab3}%
\centering
\begin{tabular}
[c]{cc}%
\toprule Parameter & Sensitivity index $\xi_{k}^{\mathcal{R}_{c}}$ of
$\mathcal{R}_{c}$\\
\midrule $\beta$ & +1\\
$b$ & +0.5863\\
$q_{2}$ & -0.3626\\
$r_{2}$ & -0.2237\\
$a$ & +0.1910\\
$q_{1}$ & -0.1701\\
$r_{1}$ & -0.0526\\
\bottomrule &
\end{tabular}
\end{table}The control reproduction number $\mathcal{R}_{c}$ is the most
sensitive to transmission rate $\beta$ and the least sensitive to recovery
rate $r_{1}$ of symptomatic infections. Additionally, $\mathcal{R}_{c}$ is
positively correlated with $a$, $\beta$, $b$, and negatively correlated with
$q_{1}$, $q_{2}$, $r_{1}$, $r_{2}$. Theoretically, the most effective measures
to control COVID-19 in India are to increase the quarantine rate $q_{2}$ of
asymptomatic infections, reduce the transmission rate $\beta$ and $b$ of
infected individuals, and enhance the cure rate $r_{2}$ of asymptomatic infections.

\subsection{ Case study in short-term COVID-19 model}

\label{Scs} 
\begin{figure}[ptbh]
\centering
\begin{minipage}[t]{0.45\linewidth}
\centering
\includegraphics[scale=0.52]{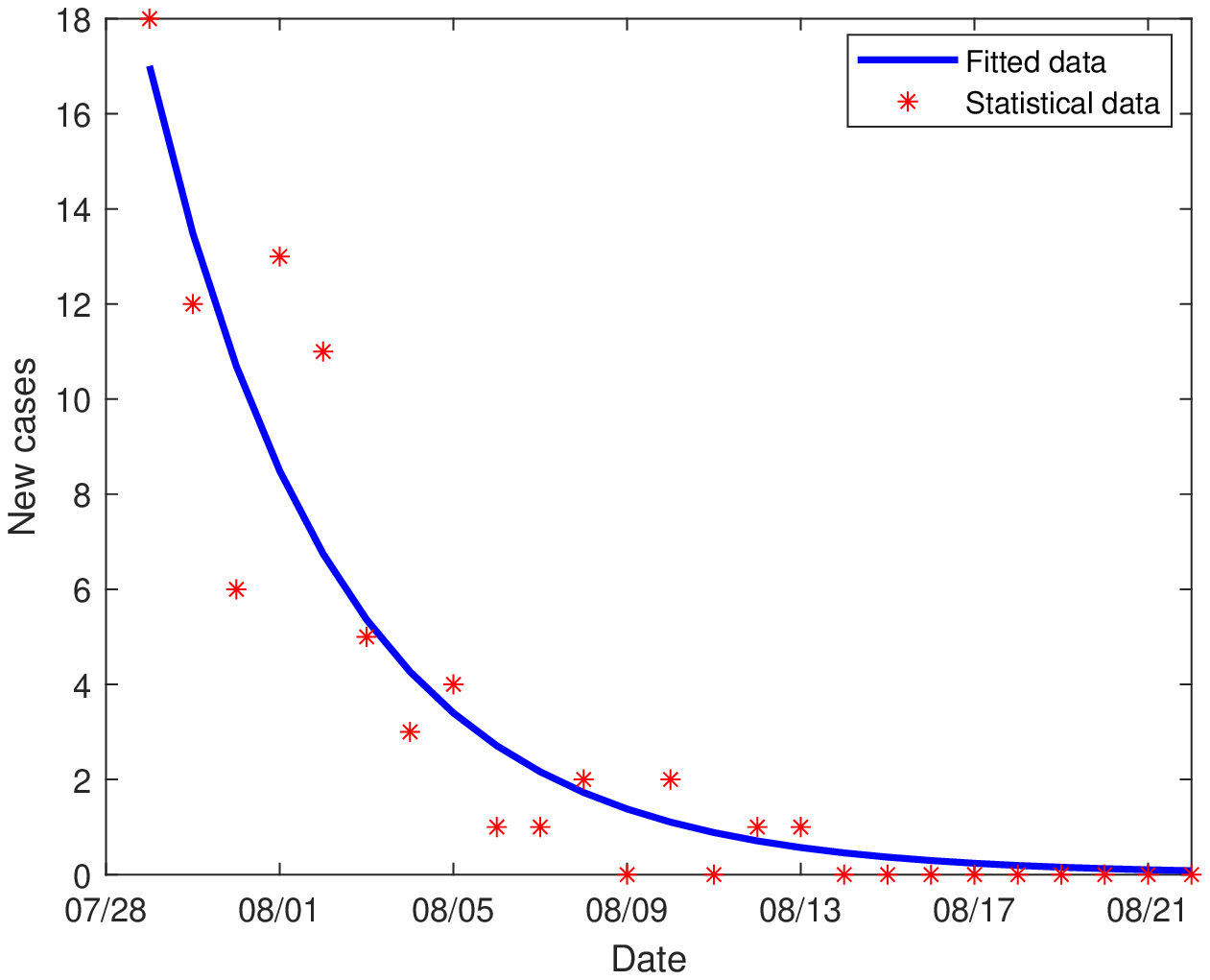}
\caption{The daily new cases in Nanjing.}\label{fig2}
\end{minipage}\begin{minipage}[t]{0.45\linewidth}
\centering
\includegraphics[scale=0.52]{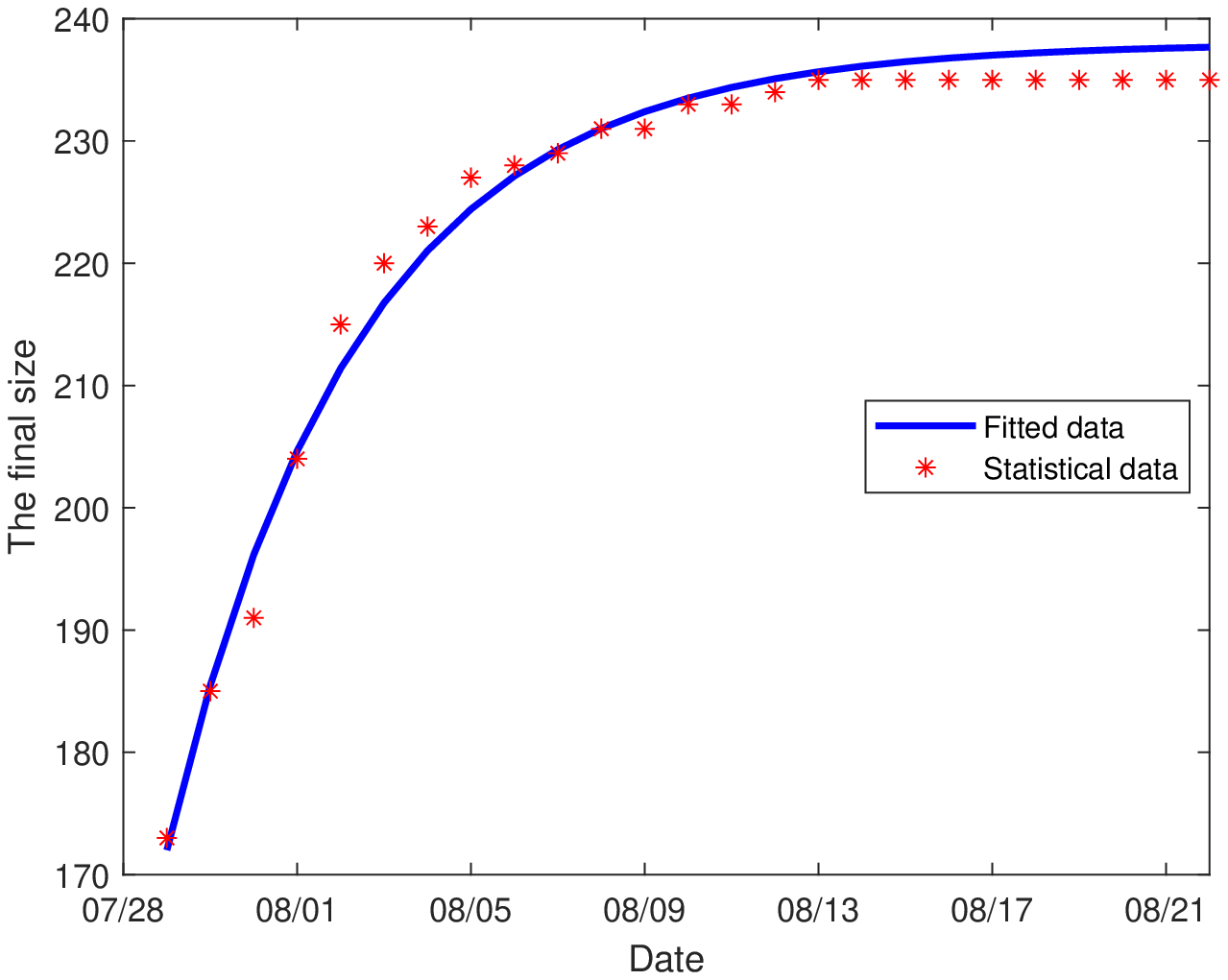}
\caption{The final size of Nanjing.}\label{fig3}
\end{minipage}
\end{figure}

Nanjing Bureau of Statistics reported \cite{Nanjing21} on 24 May
2021 that population of Nanjing is $\mathcal{N}=9314685$. We collected the
number of daily new COVID-19 cases and cumulative confirmed COVID-19 cases
from \cite{NMHC20}. Furthermore, we collated and processed the data in
\cite{NMHC20} and we considered all patients with positive nucleic acid to be
infected. Nanjing showed a sudden surge of daily new COVID-19 cases on July
25-27, 2021 since mass nucleic acid testing of Nanjing in the second round
\cite{NC25}. Consequently, to get good fitting parameter results, we selected
the daily new case data after July 27, 2021, and estimated some parameter
values of Nanjing, see Tab. \ref{tab:parameter2}.
\begin{table}[ptbh]
\caption{Parameter estimates of COVID-19 model \eqref{mod7}.}%
\label{tab:parameter2}%
\centering
\begin{tabular}
[c]{cccc}%
\toprule[0.6mm] Parameter & Parameter value of Nanjing & Source\\\hline
$\beta$ & $0.0020 $ & Estimated \\
$a $ & 0.0168 & Estimated \\
$b $ & 0.0090 & Estimated\\
$c $ & $\frac{1}{5.2}$ & \cite{Li20} \\
$p $ & 0.8702 & \cite{Gao07} \\
$q_{1} $ & 0.0757 & Estimated \\
$q_{2} $ & 0.1701 & Estimated \\
$r_{1} $ & 0.1621 & Estimated \\
$r_{2} $ & $\frac{1}{12.5} $ & \cite{Ming20} \\
$r_{3} $ & $\frac{1}{10} $ & \cite{Shao20} \\\hline\hline
$S(0) $ & $9.3144\times10^{6}$ & Estimated \\
$E(0) $ & 9.9246 & Estimated \\
$I(0) $ & $1.3604\times10^{2}$ & Estimated \\
$A(0) $ & 39.423 & Estimated \\
$Q(0) $ & 71.140 & Estimated \\
$R(0) $ & 9.7078 & Estimated \\
\bottomrule[0.6mm] &  &
\end{tabular}
\end{table}
The daily new cases and the
final size in Nanjing were further fitted (blue part) and compared with the
statistical data (red part) as shown in Figs. \ref{fig2} and \ref{fig3}. The
COVID-19 in Nanjing broke out on July 20, 2021 and disappeared on August 13,
2021. 
As shown in Fig.~\ref{fig3}, the statistical value of the final size is 235, while the ultimately
fitted value is about 238, and the relative error is about $1.28\%$,
which indicates that our model is applicable.

\begin{figure}[ptbh]
\centering
\begin{minipage}[t]{0.45\linewidth}
\centering
\includegraphics[scale=0.52]{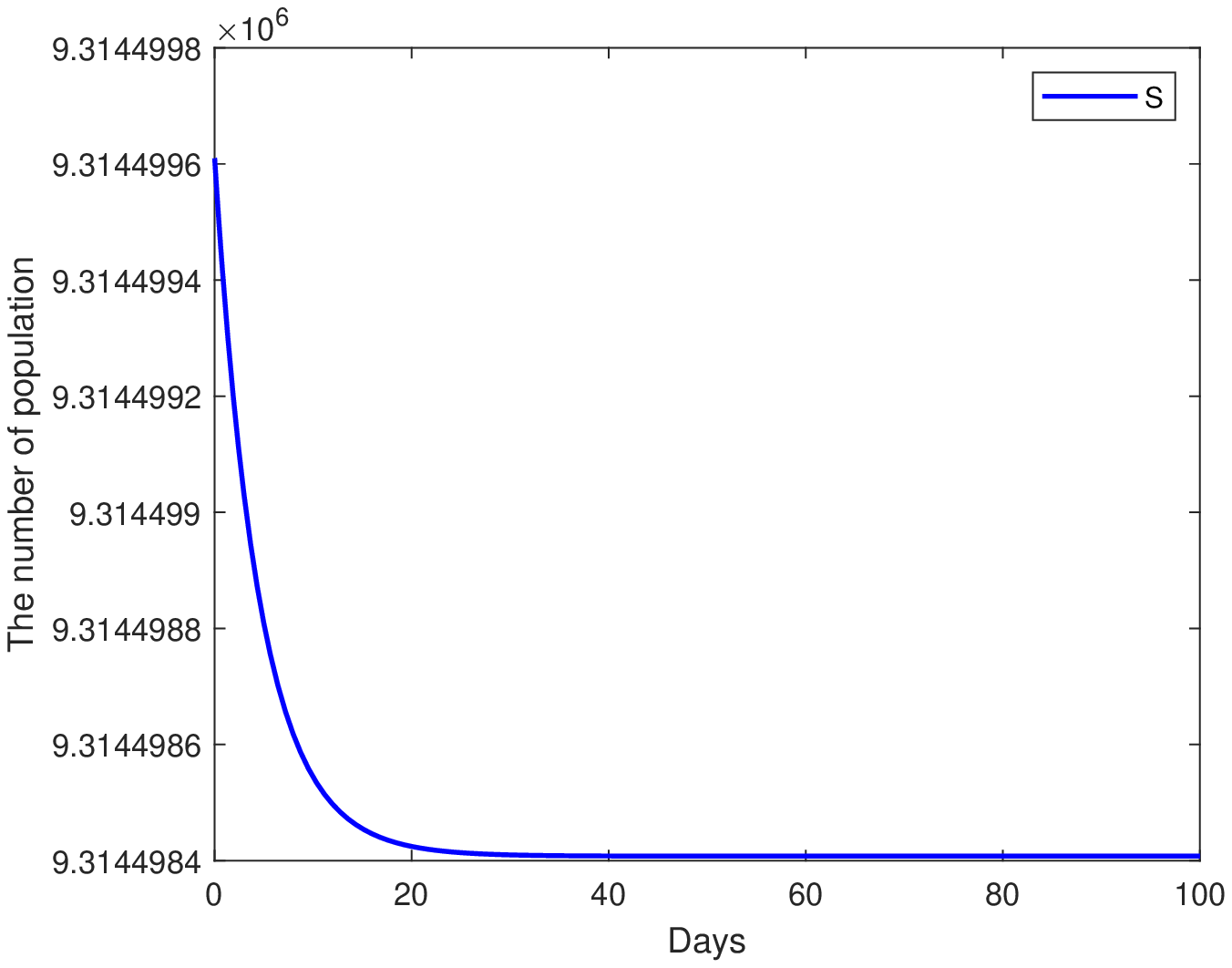}
\caption{Time evolution of population in S \\ compartment in Nanjing.}\label{fig4}
\end{minipage}\begin{minipage}[t]{0.45\linewidth}
\centering
\includegraphics[scale=0.52]{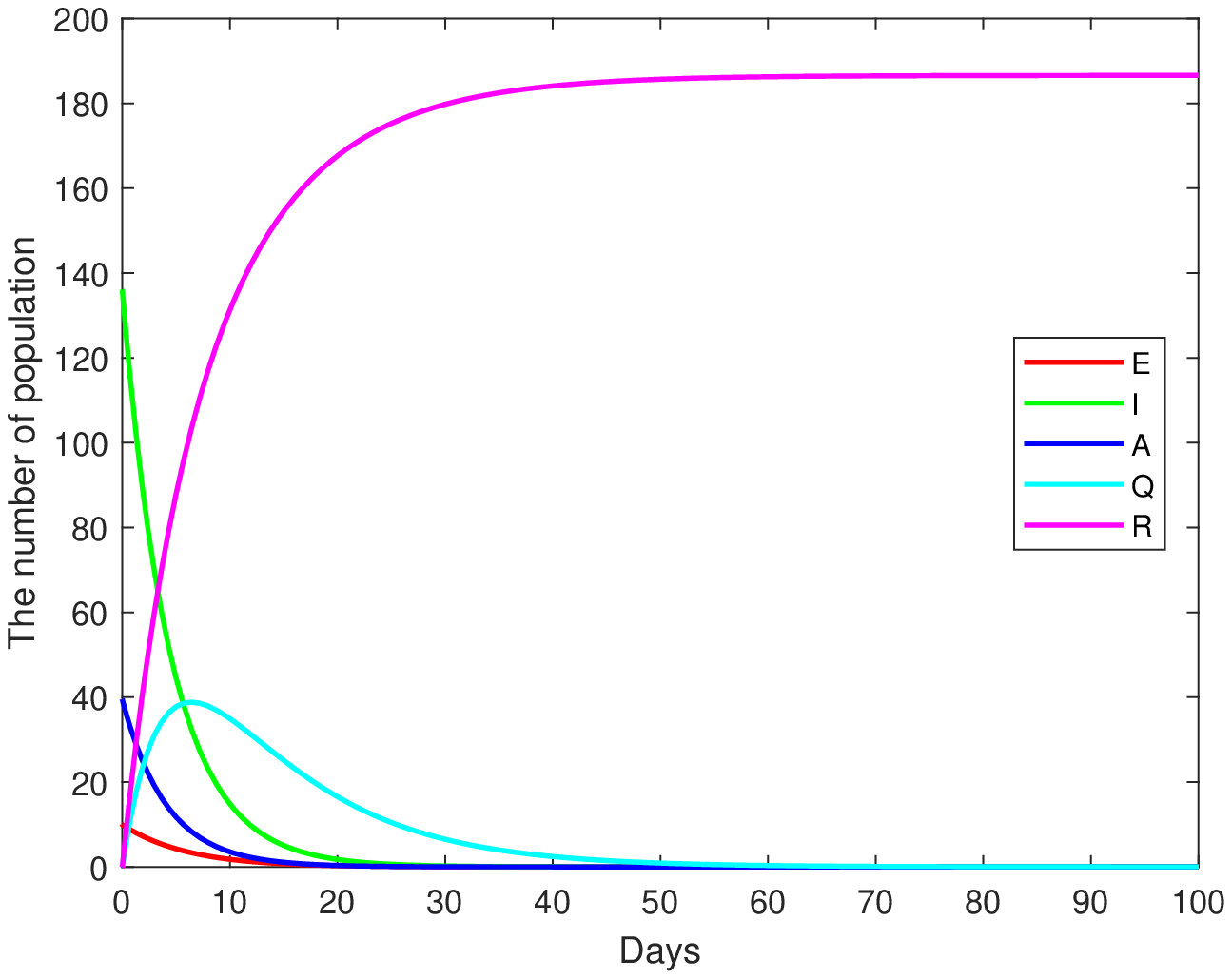}
\caption{Time evolution of population in different \\ compartments in Nanjing.}\label{fig5}
\end{minipage}
\end{figure}

The control reproduction number $\mathcal{R}_{c}$ in Nanjing can be calculated
to be $0.0073<1$. From Figs. \ref{fig4} and \ref{fig5}, we can see that
exposed individuals, symptomatically infected individuals, asymptomatically
infected individuals and quarantined individuals will eventually tend to zero,
which means that the COVID-19 will eventually disappear, see Theorems
\ref{thmS1} and \ref{thmS3}. 

\subsubsection{The impact of quarantine measures on COVID-19}

\begin{figure}[ptbh]
\centering
\includegraphics[scale=0.52]{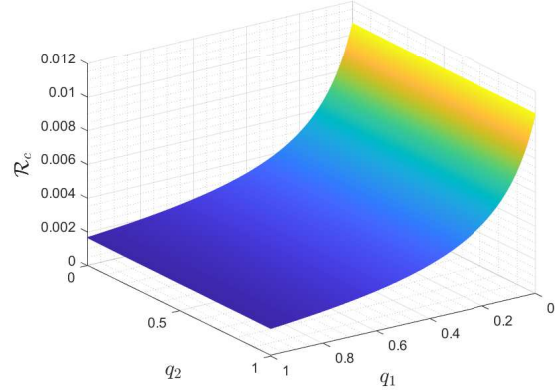} \caption{The relationship
among $\mathcal{R}_{c}$, $q_{1}$ and $q_{2}$ in Nanjing.}%
\label{fig11}%
\end{figure}

\begin{figure}[ptb]
\centering
\begin{minipage}[t]{0.45\linewidth}
\centering
\includegraphics[scale=0.52]{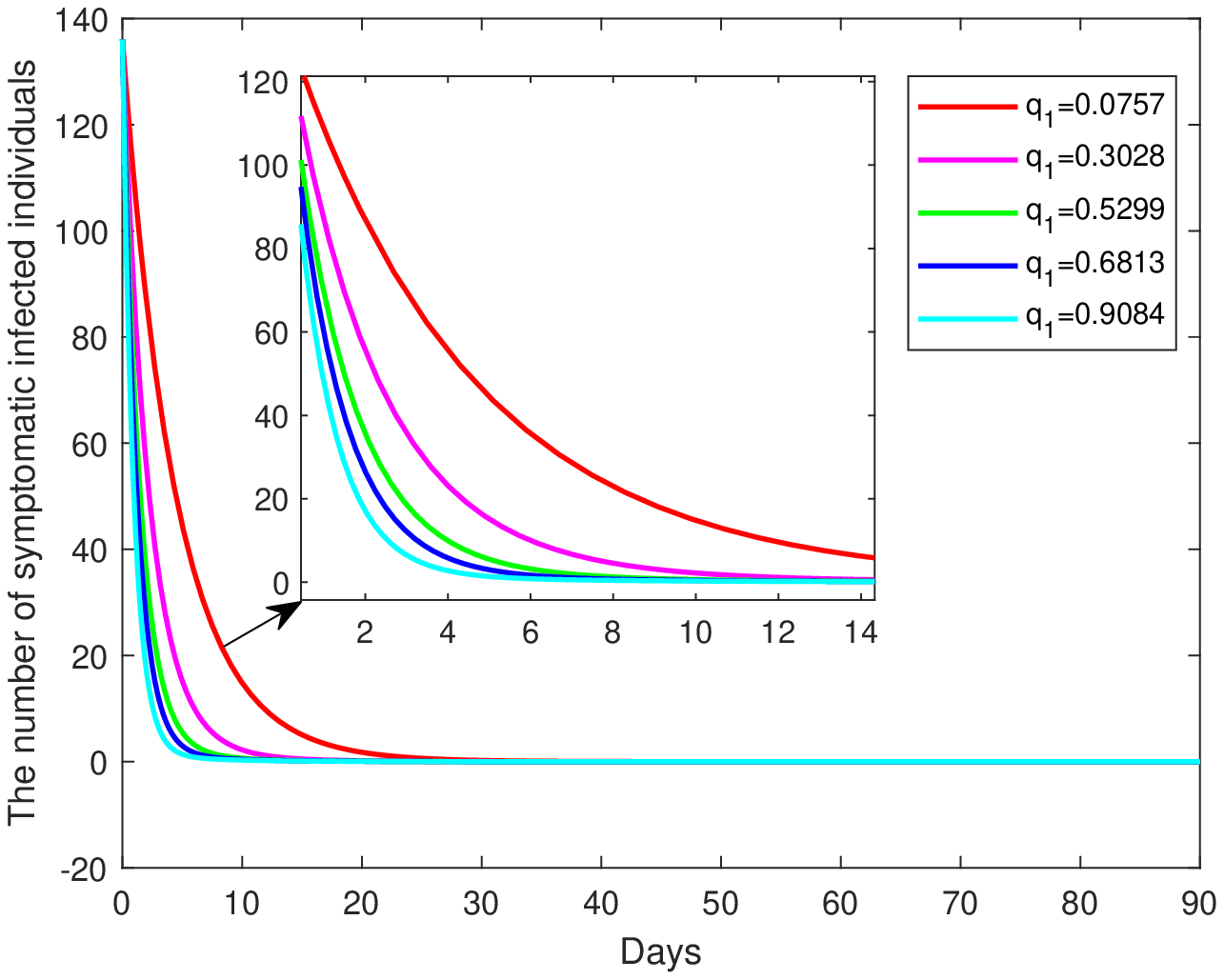}
\caption{The relationship between $I(t)$ and \\$q_{1}$ in model \eqref{mod7}.}\label{fig12}
\end{minipage}\begin{minipage}[t]{0.45\linewidth}
\centering
\includegraphics[scale=0.52]{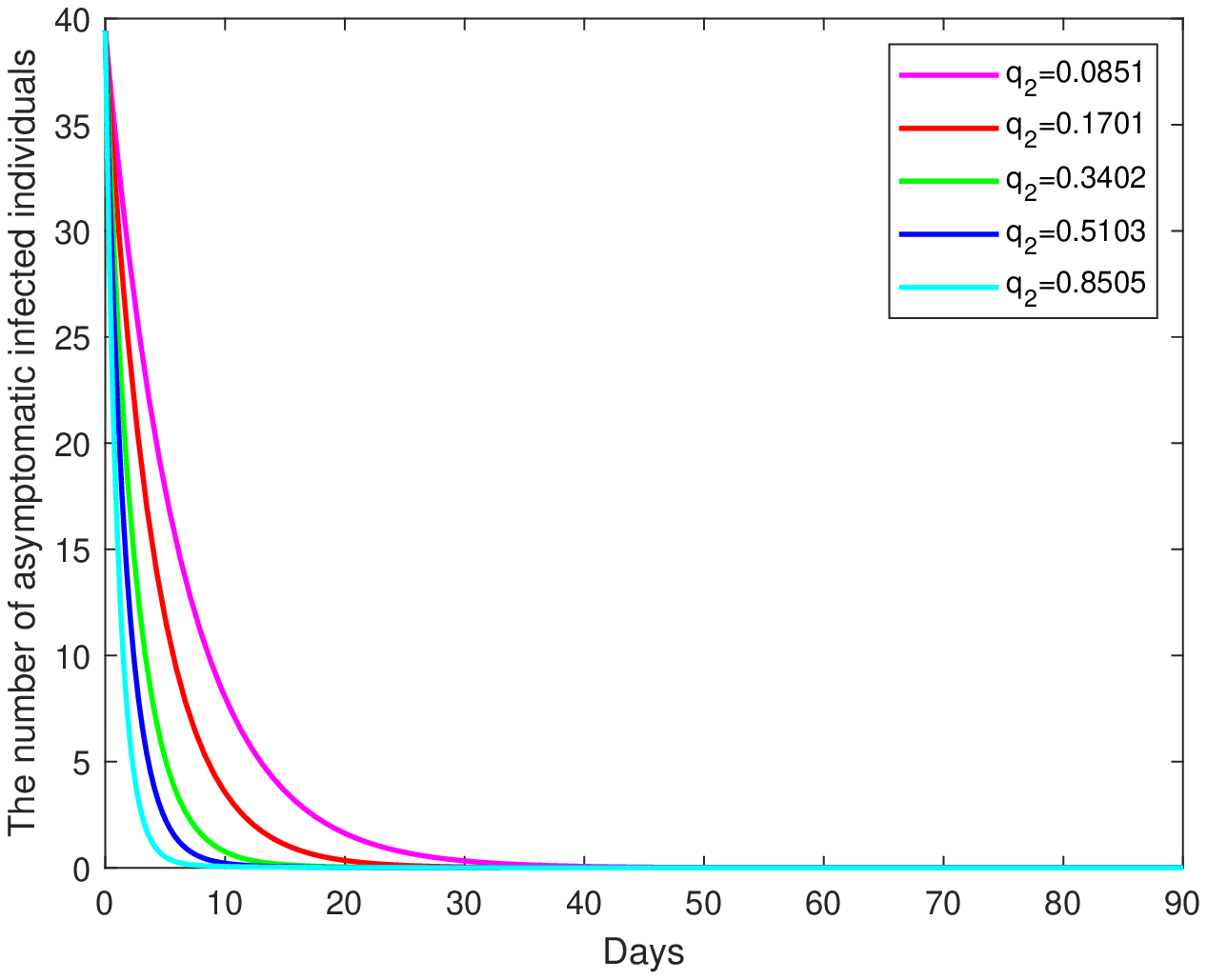}
\caption{The relationship between $A(t)$ and \\$q_{2}$ in model \eqref{mod7}.}\label{fig13}
\end{minipage}
\end{figure}

\begin{figure}[ptb]
\centering
\begin{minipage}[t]{0.45\linewidth}
\centering
\includegraphics[scale=0.52]{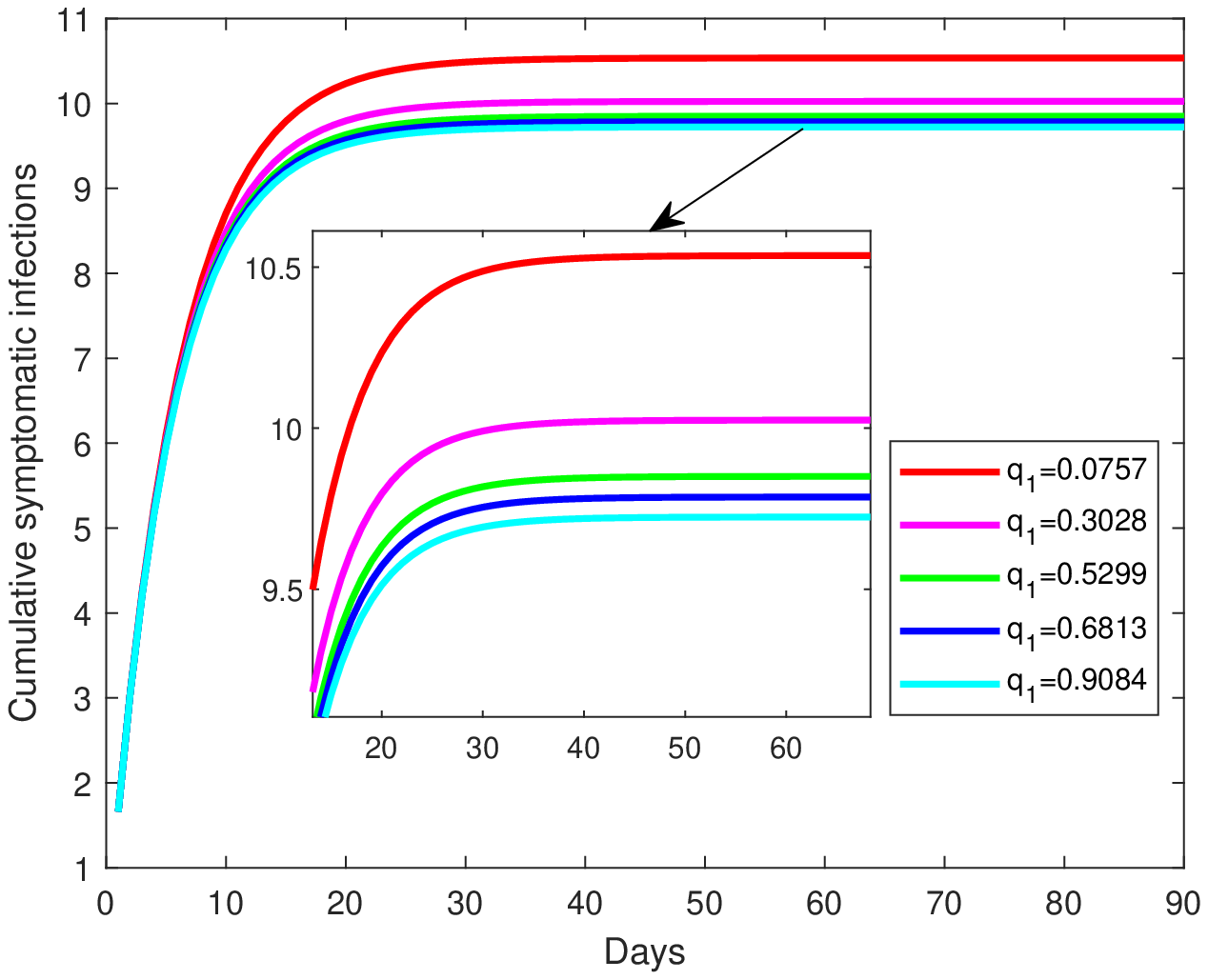}
\caption{The cumulative symptomatic \\infections in model \eqref{mod7}.}\label{fig121}
\end{minipage}\begin{minipage}[t]{0.45\linewidth}
\centering
\includegraphics[scale=0.52]{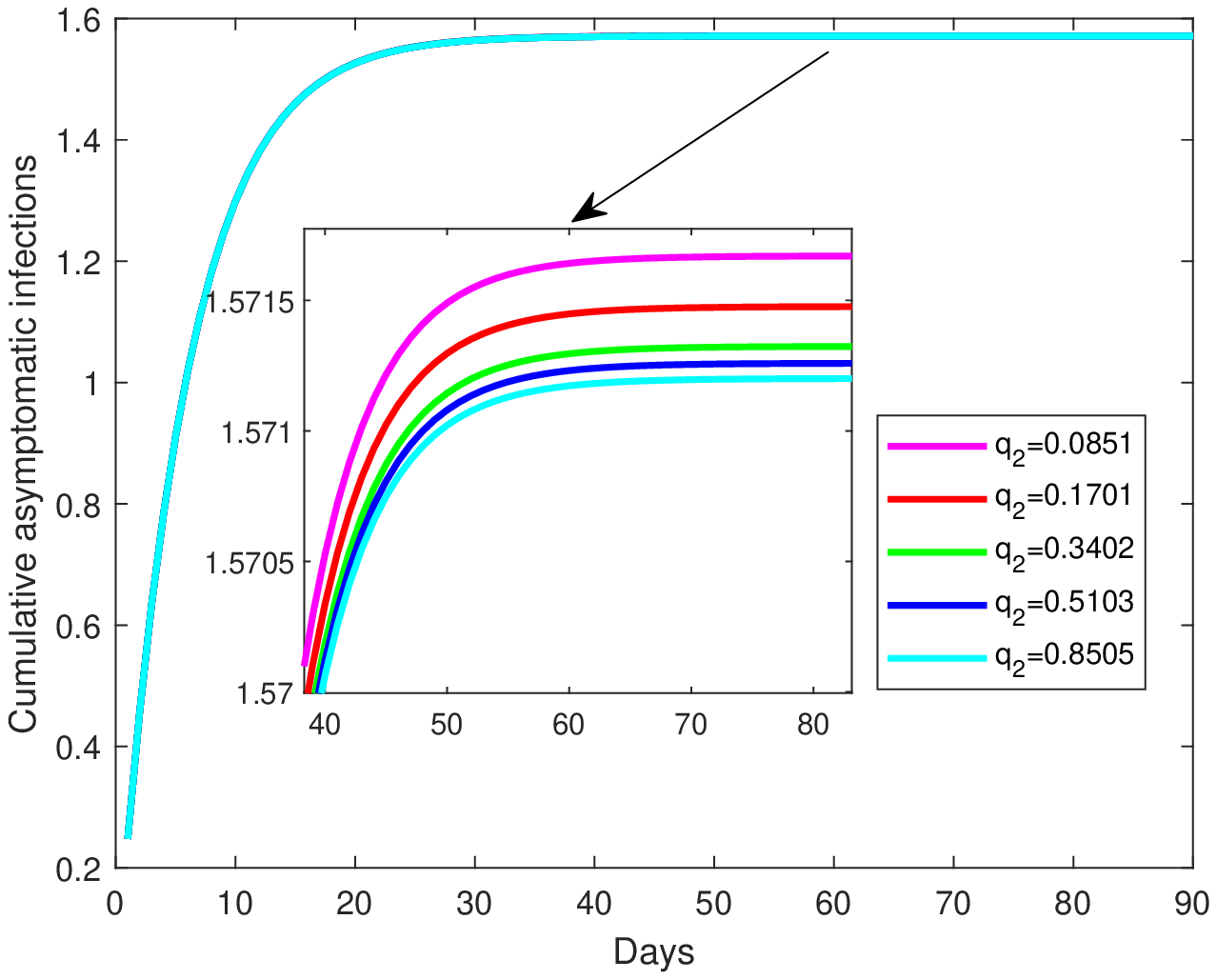}
\caption{The cumulative asymptomatic \\infections in model \eqref{mod7}.}\label{fig131}
\end{minipage}
\end{figure}

From Fig. \ref{fig11}, we find that as the values of $q_{1}$ and $q_{2}$
increase, the value of $\mathcal{R}_{c}$ decreases. And it is more sensitive
to change in $q_{1}$ than that of $q_{2}$, which is consistent with our discussion in
Section \ref{ScsSA}. Therefore, for COVID-19 in Nanjing at that time,
strengthening quarantine measure of symptomatically infected individuals is
more necessary. In Figs. \ref{fig12} and \ref{fig121}, the value of $q_{1}$ is
changed while other parameters are fixed. We find that with the increase in $q_{1}$, $I(t)$ tends to zero faster, as shown in Fig. \ref{fig12},
and the cumulative confirmed symptomatic infections also decrease (see
Fig. \ref{fig121}). In Figs. \ref{fig13} and \ref{fig131}, only the value of
$q_{2}$ changes. As the value of $q_{2}$ increases, the time when $A(t)$
approaches zero is shorter and the cumulative
number of asymptomatic infections shows a downward trend, as shown in Fig. \ref{fig13} and Fig.
\ref{fig131}, respectively. This suggests that the stronger the quarantine measures, the
less harm COVID-19 will cause to human beings. 

\subsubsection{The impact of asymptomatic infections on COVID-19}

\label{622}
\begin{figure}[ptb]
\centering
\begin{minipage}[t]{0.45\linewidth}
\centering
\includegraphics[scale=0.52]{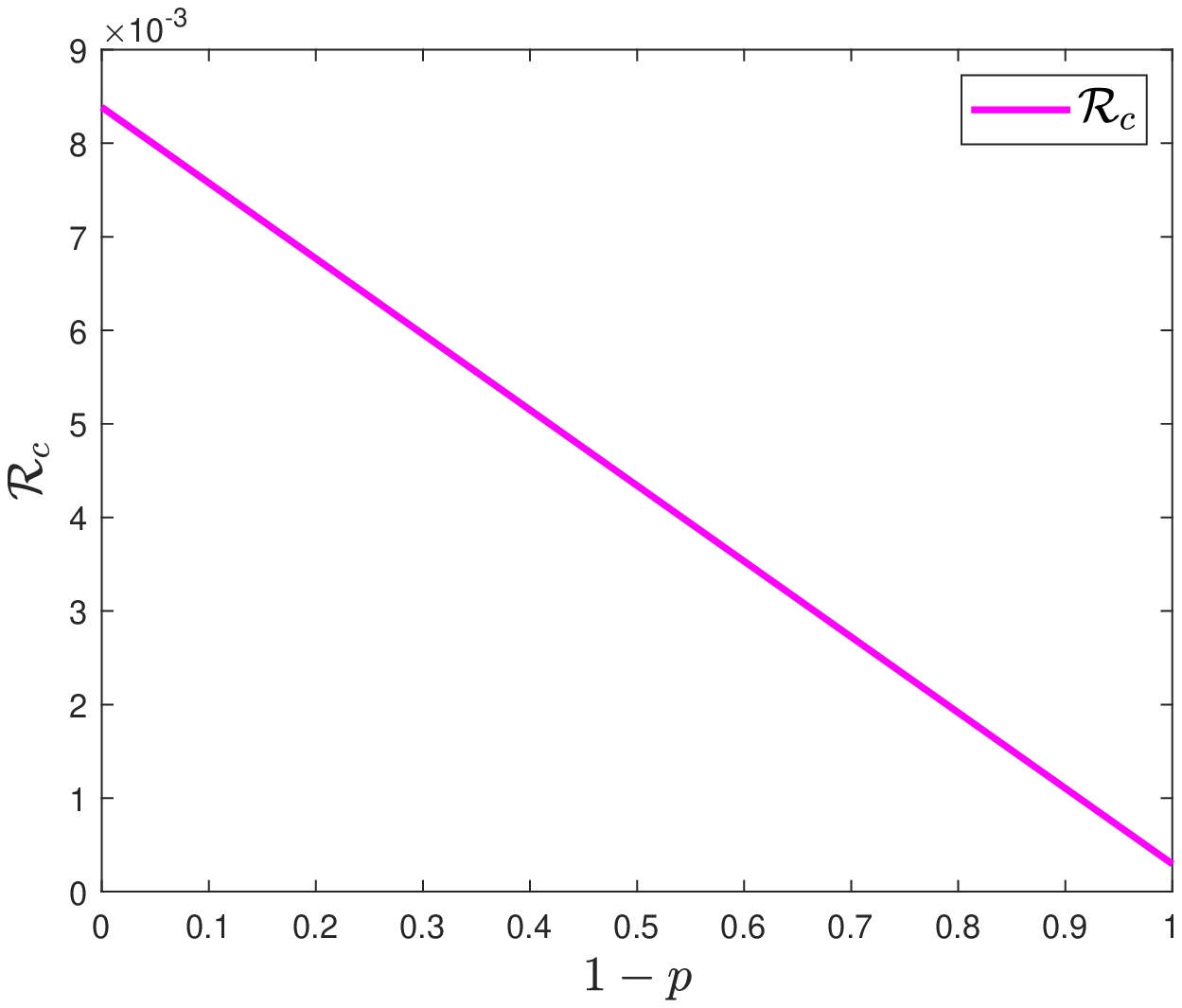}
\caption{Relationship between $\mathcal{R}_{c}$ and $1-p$.}\label{figSA1}
\end{minipage}\begin{minipage}[t]{0.45\linewidth}
\centering
\includegraphics[scale=0.52]{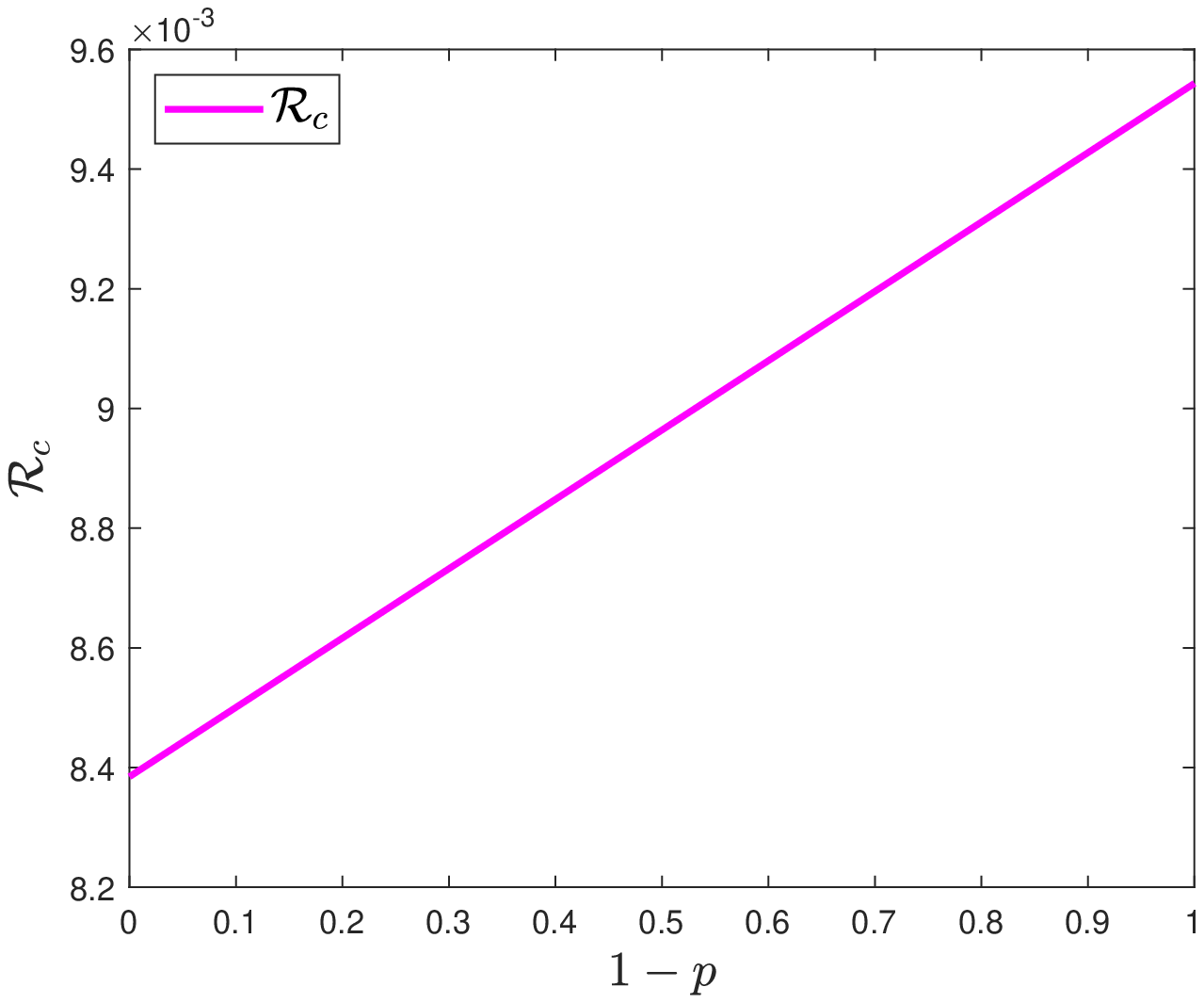}
\caption{Relationship between $\mathcal{R}_{c}$ and $1-p$.}\label{figSA2}
\end{minipage}
\end{figure}

\begin{figure}[ptb]
\centering
\includegraphics[scale=0.52]{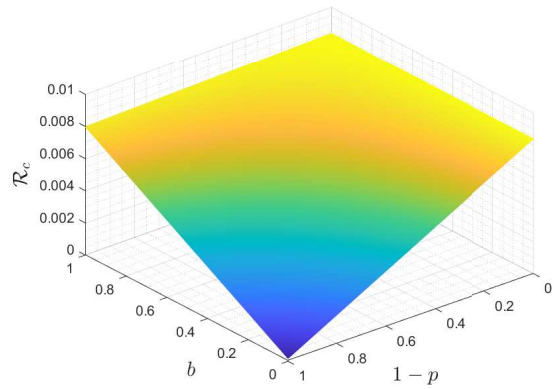} \caption{The
relationship among $\mathcal{R}_{c}$, $b$ and $1-p$ in Nanjing.}%
\label{figSA3}%
\end{figure}

In the outbreak of COVID-19 in Nanjing, asymptomatic infections
accounted for a small proportion of the total infections \cite{NMHC20}.
Nonetheless, asymptomatic transmission had played a non-negligible role in that outbreak.
According to the data in Tab. \ref{tab:parameter2}, it can be calculated that
the critical value of $b$ is about 1.0517. Hence, when $b\leq1.0517$,
$\mathcal{R}_{c}$ is negatively correlated with $1-p$, as shown in Fig.
\ref{figSA1}; when $b>1.0517$, $\mathcal{R}_{c}$ is positively correlated with
$1-p$, as shown in Fig. \ref{figSA2}. Meanwhile, Fig. \ref{figSA3} visually
shows the effect of changes in $1-p$ and $b$ on $\mathcal{R}_{c}$.


\subsubsection{Sensitivity analysis}

\label{ScsSA} From \eqref{scrn} and \eqref{tsi}, we can obtain the sensitivity
index expressions of $\mathcal{R}_{c}$ with respect to seven parameters,
respectively. The sensitivity index can be obtained by using the data in Tab.
\ref{tab:parameter2}, as shown in Tab. \ref{tab4}. The parameters in the table
are arranged from the most sensitive to the least sensitive.
\begin{align*}
\xi_{a}^{\mathcal{R}_{c}}  &  =\frac{\beta}{c}\cdot\frac{a}{\mathcal{R}_{c}}
=\frac{a(q_{1}+r_{1})(q_{2}+r_{2})}{a(q_{1}+r_{1})(q_{2}+r_{2})+pc(q_{2}%
+r_{2})+bc(1-p)(q_{1}+r_{1})},\\
\xi_{\beta}^{\mathcal{R}_{c}}  &  =\left(  \frac{a}{c}+\frac{p}{q_{1}+r_{1}}
+\frac{b(1-p)}{q_{2}+r_{2}}\right)  \cdot\frac{\beta}{\mathcal{R}_{c}}=1,\\
\xi_{b}^{\mathcal{R}_{c}}  &  =\frac{\beta(1-p)}{q_{2}+r_{2}} \cdot\frac
{b}{\mathcal{R}_{c}} =\frac{cb(1-p)(q_{1}+r_{1})} {a(q_{1}+r_{1})(q_{2}%
+r_{2})+pc(q_{2}+r_{2})+cb(1-p)(q_{1}+r_{1})},\\
\xi_{q_{1}}^{\mathcal{R}_{c}}  &  =-\frac{p\beta}{(q_{1}+r_{1})^{2}}\cdot
\frac{q_{1}}{\mathcal{R}_{c}} =-\frac{pcq_{1}(q_{2}+r_{2})}{(q_{1}%
+r_{1})\left[  a(q_{1}+r_{1})(q_{2}+r_{2})+pc(q_{2}+r_{2})+cb(1-p)(q_{1}%
+r_{1})\right]  },\\
\xi_{q_{2}}^{\mathcal{R}_{c}}  &  =-\frac{b\beta(1-p)}{(q_{2}+r_{2})^{2}}%
\cdot\frac{q_{2}}{\mathcal{R}_{c}} =-\frac{bc(1-p)q_{2}(q_{1}+r_{1})}%
{(q_{2}+r_{2})\left[  a(q_{1}+r_{1})(q_{2}+r_{2})+pc(q_{2}+r_{2}%
)+cb(1-p)(q_{1}+r_{1})\right]  },\\
\xi_{r_{1}}^{\mathcal{R}_{c}}  &  =-\frac{p\beta}{(q_{1}+r_{1})^{2}}\cdot
\frac{r_{1}}{\mathcal{R}_{c}} =-\frac{pcr_{1}(q_{2}+r_{2})}{(q_{1}%
+r_{1})\left[  a(q_{1}+r_{1})(q_{2}+r_{2})+pc(q_{2}+r_{2})+cb(1-p)(q_{1}%
+r_{1})\right]  },\\
\xi_{r_{2}}^{\mathcal{R}_{c}}  &  =-\frac{b\beta(1-p)}{(q_{2}+r_{2})^{2}}%
\cdot\frac{r_{2}}{\mathcal{R}_{c}} =-\frac{bc(1-p)r_{2}(q_{1}+r_{1})}%
{(q_{2}+r_{2})\left[  a(q_{1}+r_{1})(q_{2}+r_{2})+pc(q_{2}+r_{2}%
)+cb(1-p)(q_{1}+r_{1})\right]  }.
\end{align*}
\begin{table}[ptbh]
\caption{Sensitivity index of $\mathcal{R}_{c}$ with respect to parameters.}%
\label{tab4}%
\centering
\begin{tabular}
[c]{ccc}%
\toprule Parameter & Sensitivity index $\xi_{k}^{\mathcal{R}_{c}}$ of
short-term model $\mathcal{R}_{c}$ & \\
\midrule $\beta$ & +1 & \\
$r_{1}$ & -0.6648 & \\
$q_{1}$ & -0.3106 & \\
$a$ & +0.0233 & \\
$b$ & +0.0012 & \\
$q_{2}$ & -0.0008 & \\
$r_{2}$ & -0.0004 & \\
\bottomrule &  &
\end{tabular}
\end{table}Similar to the conclusion in Tab. \ref{tab3}, the spread of the
pandemic can be controlled by reducing transmission rates $a$, $\beta$, $b$,
strengthening quarantine measures $q_{1}$, $q_{2}$, and enhancing recovery
rates $r_{1}$, $r_{2}$. The most effective measures for Nanjing were to reduce
the transmission of symptomatically infected individuals and enhance the
recovery rate of symptomatic infections and quarantine measures.

\section{Conclusions}

This paper not only analyzes the global stability of the COVID-19-free
equilibrium $V^{0}$ and the COVID-19 equilibrium $V^{\ast}$ of model
\eqref{mod1}, but also solves the left problems in \cite{Guo22}. For the local
stability of the COVID-19 equilibrium $V^{\ast}$, it is difficult to use the
Routh-Hurwitz criterion. To this end, we make use of proof by contradiction
and the properties of complex modulus with some novel techniques and less
computation. It is well known that the persistence result of model \eqref{mod1} is
essential for the global attractivity of $V^{\ast}$, and hence we prove weak
persistence of model \eqref{mod1}.
To obtain the global stability results of $V^{0}$ and $V^{\ast}$,
we adopt the limit system of model \eqref{mod1} and Lyapunov function method.
Specifically, $V^{0}$ is globally asymptotically stable for $\mathcal{R}%
_{c}<1$ and globally attractive for $\mathcal{R}_{c}=1$ in $\mathbb{R}_{+}%
^{6}$, which implies that COVID-19 will disappear;
$V^{\ast}$ is globally asymptotically stable for $\mathcal{R}_{c}>1$ in $\Omega$, which indicates
that COVID-19 will persist.

Although COVID-19 is in long-term development around the world since December
2019, COVID-19 disappeared in a short term on account of the
strong prevention and control measures in some areas. For this reason, we propose a
short-term COVID-19 model \eqref{mod7} based on model \eqref{mod1}. It is sure
that model \eqref{mod7} has no COVID-19 equilibrium. We work out the stability
of the multiple COVID-19-free equilibria and the expression of the final size.

A retrospective study was conducted on the transmission of COVID-19 in India and Nanjing.
We apply the long-term and the short-term models with publicly
available official statistics to numerically demonstrate different transmission
characteristics of COVID-19 in India and Nanjing, respectively. These numerical simulations validate the theoretical results of Theorems \ref{thm51}, \ref{thm52}, \ref{thmS1} and \ref{thmS3}. The long-term model can well predict the spread of COVID-19 in India, and the short-term model perfectly fits the final size of Nanjing. Particularly, the relationship between the proportion of asymptomatic infected individuals and $\mathcal{R}_{c}$ is related to the transmission ability of asymptomatic infected individuals. As for India,
case study shows that enhanced quarantine measures can not only prevent COVID-19
transmission but also reduce the peak value of infected individuals and
cumulative confirmed cases. Asymptomatic transmission played a key role in
the outbreak in India. For Nanjing, compared with asymptomatic infections,
strengthening the quarantine measures for symptomatic infections is more
conducive to controlling the spread of the COVID-19 at that time. With the
strengthening of quarantine measures, infected individuals will tend to zero
faster and the final COVID-19 size will decrease. It can be found that, although
there are a few asymptomatic infections of the COVID-19 in Nanjing, the impact
of asymptomatic transmission cannot be ignored, otherwise it will lead to
inaccurate calculation of control reproduction number in Nanjing. As shown in Tabs. \ref{tab3} and \ref{tab4},
the spread of COVID-19 can be reduced by vaccination, wearing masks, mass
nucleic acid testing and strengthening quarantine measures, and so on. In
addition, some measures to speed up the recovery rate of patients are also
beneficial to control the COVID-19 pandemic.

\section*{Acknowledgements}

This work is partially supported by the National NSF of China (Nos. 11901027,
11871093, 11671382, 11971273 and 12126426), the Major Program of the National
NSF of China (No. 12090014), the State Key Program of the National NSF of
China (No. 12031020), the NSF of Shandong Province (No. ZR2018MA004), and the
China Postdoctoral Science Foundation (No. 2021M703426), the Pyramid Talent
Training Project of BUCEA (No. JDYC20200327), and the BUCEA Post Graduate
Innovation Project (No. PG2022143). The authors would like to thank Prof.
Jing-An Cui and Dr. Liping Sun for their helpful suggestions.

\section*{References}


\end{document}